\documentclass[12pt]{amsart}
\usepackage{amssymb}
\usepackage{amsmath}
\usepackage{pictex}
\def\func#1{\mathop{\rm #1}}%

\newcommand\A{\mathbb A}
\newcommand\Af{\mathcal{A}}
\newcommand\Aff{\operatorname{\sf Aff}}
\newcommand\Alt{\mathfrak{A}}
\newcommand\AND{\quad\text{and}\quad}
\newcommand\asf{\mathsf{a}}
\newcommand\Aut{\operatorname{\sf Aut}}
\newcommand\bd{\partial}
\newcommand\bsf{\mathsf{b}}
\newcommand\C{\mathbb C}
\newcommand\cf{\curlywedge}
\newcommand\de{\delta}
\newcommand\DL{\mathsf{DL}}
\newcommand\dps{\displaystyle}
\newcommand\ee{\mathbf e}

\newcommand\F{\mathbb F}
\newcommand\ff{\text{\sl f\/}}
\newcommand\fsf{\mathsf{f}}

\newcommand\Ga{\Gamma}
\newcommand\geo[1]{\overline{#1}}

\newcommand\gsf{\mathsf{g}}

\newcommand\he{N}

\newcommand\hh{\mathbf h}
\newcommand\hor{\mathfrak{h}}
\newcommand\Hor{\mathfrak{H}}
\newcommand\id{\text{\sl id}}
\newcommand\im{\mathfrak{i}\,}
\newcommand\Int{\mathfrak O}

\newcommand\kk{\mathbf k}
\newcommand\Lb{\mathfrak L}
\newcommand\Lq{\Lb_q}
\newcommand\la{\lambda}
\newcommand\lk{\mathbf l}
\newcommand\lle{\preccurlyeq}
\newcommand\MC{\mathcal{M}}
\newcommand\mm{\mathbf m}
\newcommand\N{\mathbb N}
\newcommand\Oct{{\mathcal O}}
\newcommand\om{\omega}
\newcommand\Pol{\mathsf{S}}
\newcommand\Prob{\mathsf{Pr}}
\newcommand\Psf{\mathsf{P}}
\newcommand\Q{\mathbb{Q}}
\newcommand\Qsf{\mathsf{Q}}
\newcommand\R{\mathbb R}
\newcommand\Ring{\mathcal R}
\newcommand\SF{\mathfrak{S}}

\newcommand\si{\sigma}
\newcommand\spec{\operatorname{\sf spec}}
\newcommand\spn{\operatorname{\sf span}}

\newcommand\SUM{\mathsf{Z}}

\newcommand\tsf{\mathsf{t}}
\newcommand\T{\mathbb T}
\newcommand\tb{\mathbf{t}}
\newcommand\up{\mathfrak{u}}
\newcommand\uno{\mathbf{1}}
\newcommand\Vf{\mathfrak{V}}
\newcommand\Vol{\operatorname{\sf Vol}}
\newcommand\vsf{\mathsf{v}}
\newcommand\wh{\widehat}
\newcommand\wt{\widetilde}
\newcommand\Z{\mathbb Z}
\newcommand\zero{\mathbf{0}}
\newcommand\Zyk{\mathfrak Z}
\newcommand\Zq{\Zyk_q}

\numberwithin{equation}{section}

\newtheoremstyle{mythm}
  {9pt}
  {9pt}
  {\itshape}
  {0pt}
  {\bfseries}
  {}
  { }
  {\thmnumber{(#2)}\thmname{ #1}\thmnote{ #3}}

\newtheoremstyle{mydef}
  {9pt}
  {9pt}
  {\normalfont}
  {0pt}
  {\bfseries}
  {}
  { }
  {\thmnumber{(#2)}\thmname{ #1}\thmnote{ #3}}

\theoremstyle{mythm}
\newtheorem{thm}[equation]{Theorem.}
\newtheorem{pro}[equation]{Proposition.}
\newtheorem{lem}[equation]{Lemma.}
\newtheorem{cor}[equation]{Corollary.}

\theoremstyle{mydef}
\newtheorem{dfn}[equation]{Definition.}

\newtheorem{rmk}[equation]{Remark.}

\begin{document}$\,$ \vspace{-1truecm}
\title{\large  Horocyclic products of trees}
\author{\bf Laurent BARTHOLDI, Markus NEUHAUSER and Wolfgang WOESS}
\address{\parbox{.8\linewidth}{\'Ecole Polytechnique F\'ed\'erale de 
    Lausanne (EPFL),\\ Institut de Math\'ematiques B (IMB),\\
    CH-1015 Lausanne, Switzerland\\}}
\email{laurent.bartholdi@epfl.ch}
\address{\parbox{.8\linewidth}{Mathematisches Institut,
Georg-August-Universit\"{a}t G\"{o}ttingen,\\
Bunsenstra\ss e 3--5,
D-37073 G\"{o}ttingen, 
Germany\\}}
\email{neuhause@uni-math.gwdg.de}
\address{\parbox{.8\linewidth}{Institut f\"ur Mathematik C, 
Technische Universit\"at Graz,\\
Steyrergasse 30, A-8010 Graz, Austria\\}}
\email{woess@weyl.math.tu-graz.ac.at}
\date{17 January 2006} 
\thanks{Supported by FWF (Austrian Science Fund) project P15577}
\subjclass[2000] {05C50, 20E22, 47A10, 60B15}
\keywords{restricted wreath product, trees, horocycles, 
Diestel-Leader graph, growth function,
normal form, Markov operator, spectrum}
\begin{abstract}
  Let $T_1,\dots, T_d$ be homogeneous trees with degrees $q_1+1,
  \dots, q_d+1 \ge 3,$ respectively.  For each tree, let $\hor:T_j \to
  \Z$ be the Busemann function with respect to a fixed boundary point
  (end).  Its level sets are the horocycles. The horocyclic product of
  $T_1,\dots, T_d$ is the graph $\DL(q_1,\dots,q_d)$ consisting of all
  $d$-tuples $x_1 \cdots x_d \in T_1 \times \dots \times T_d$ with
  $\hor(x_1)+\dots+\hor(x_d)=0$, equipped with a natural neighbourhood
  relation.  In the present paper, we explore the geometric,
  algebraic, analytic and probabilistic properties of these graphs and
  their isometry groups. If $d=2$ and $q_1=q_2=q$ then we obtain a
  Cayley graph of the lamplighter group (wreath product) $\Zq \wr \Z$.
  If $d = 3$ and $q_1 = q_2 = q_3 = q$ then $\DL$ is the Cayley graph
  of a finitely presented group into which the lamplighter group
  embeds naturally. Also when $d\ge 4$ and $q_1 = \dots = q_d = q$ is
  such that each prime power in the decomposition of $q$ is larger
  than $d-1$, we show that $\DL$ is a Cayley graph of a finitely
  presented group.  This group is of type $F_{d-1}$, but not $F_d$. It
  is not automatic, but it is an automata group in most cases.  On the
  other hand, when the $q_j$ do not all coincide, $\DL(q_1,\dots,q_d)$
  is a vertex-transitive graph, but is not the Cayley graph of a
  finitely generated group. Indeed, it does not even admit a group
  action with finitely many orbits and finite point stabilizers. The
  $\ell^2$-spectrum of the ``simple random walk'' operator on $\DL$ is
  always pure point. When $d=2$, it is known explicitly from previous
  work, while for $d=3$ we compute it explicitly.  Finally, we
  determine the Poisson boundary of a large class of group-invariant
  random walks on $\DL$. It coincides with a part of the geometric
  boundary of $\DL$.
\tableofcontents
\end{abstract}

\maketitle

\markboth{{\sf L. Bartholdi, M. Neuhauser and W. Woess}}{{\sf Horocyclic products
of trees}}
\baselineskip 15pt

$\,$ \vspace{-1.5truecm}

\section{Introduction}\label{intro}
Let $X$ be a locally finite, infinite, connected graph without loops.
We write $x \sim y$ if $x,y$ are neighbours (connected by an edge),
and $\deg(x)$ for the number of neighbours of $x$.  Suppose that $X$
is written as a disjoint union of non-empty sets $H_k$, $k \in \Z$
(the \emph{horocycles}), where each element in $H_k$ has neighbours
both in $H_{k-1}$ and in $H_{k+1}$, but none in any other $H_l$. (This
condition is tailored to our purposes and can be generalized.) The
associated surjection $\hor: X \to \Z$, where $\hor(x) = k$ if $x \in
H_k$, is a graph homomorphism of $X$ onto the two-way-infinite path
$\Z$. We call it a \emph{Busemann function,} although this terminology
is justified completely only in specific cases (see below), and we say
that $(X,\hor)$ is a \emph{Busemann pair.}  Now let $X_1, \dots, X_d$
be a family of such graphs with associated Busemann functions $\hor:
X_j \to \Z$ (we use the same symbol $\hor$ for each of them). Then
their \emph{horocyclic product} is
\begin{equation}\label{eq:horoprod1}
\prod_{j=1}^d{\!_{\dps \hor}}\,\, X_j = 
\{x_1\cdots x_d \in X_1 \times \cdots \times X_d :
\hor(x_1)+ \dots +\hor(x_d) = 0 \}
\end{equation}
with neighbourhood
\begin{equation}\label{nbhd}
\begin{gathered}x_1\cdots x_d \sim y_1\cdots y_d \iff\\
\text{there are}\;\ i\neq j\;\ \text{such that}\;\ x_i \sim y_i\,,\;
x_j \sim y_j\;\ \text{and}\;\ x_k=y_k\; \text{ for all }k \ne i, j\,.
\end{gathered}
\end{equation} 
In particular, one must have $\hor(x_i)-\hor(y_i) = \hor(y_j)-\hor(x_j) 
= \pm 1$. Thus, $x=x_1\cdots x_d \mapsto \Hor(x)= \bigl(\hor(x_1),\dots,\hor(x_d)\bigr)$
is a graph homomorphism of $\prod{\!_{\dps \hor}}\,\, X_j$ onto the
simplicial lattice 
$\A_{d-1} = \{ \kk =(k_1,\dots,k_d) \in \Z^d : k_1+\dots+k_d=0\}$.
In the latter, two points are neighbours if they differ by a vector
$\ee_i-\ee_j$, where $i \ne j$ and $\ee_i \in \Z^d$ is the unit vector with a 
$1$ in its $i$-th coordinate.

There is an analogous construction for groups, compare with
{\sc Kaimanovich and Woess~\cite[p.356]{KaWo}}. Let $\Ga_1, \dots,
\Ga_d$ be topological (e.g.\, in particular, discrete)
groups, each one equipped with a
continuous homomorphism $\hor: \Ga_j \to \Z$ (or $\to \R$; we again
use the same symbol $\hor$ for each of them). Then their horocyclic
product is
\begin{equation}\label{eq:horoprod2}
\prod_{j=1}^d{\!_{\dps \hor}}\,\, \Ga_j = 
\{g_1\cdots g_d \in \Ga_1 \times \cdots \times \Ga_d :
\hor(g_1)+ \dots +\hor(g_d) = 0 \}\,,
\end{equation}
which is a closed subgroup of the direct product of the $\Ga_j$. For
finitely generated groups $\Ga_j$, this kind of construction was used
for example by {\sc Bridson}~\cite{Br}.  However, our approach has a
different ``history'', and below, the groups will be non-discrete
isometry groups of homogeneous trees.  Here, horocyclic products of
groups will arise as isometry (automorphism) groups of horocyclic
products of graphs.

If $(X_1,\hor)$ and $(X_2,\hor)$ are two Busemann pairs, then a
\emph{Busemann isometry} $g$ from the first to the latter is a graph isomorphism
$g: X_1 \to X_2$ such that $x_1 \mapsto \hor(gx_1) - \hor(x_1)$ is constant.
We write $\hor(g)$ for this constant. In particular, the group $\Aut(X,\hor)$
of a given Busemann pair $(X,\hor)$ consists of all Busemann isometries 
$X \to X$. Given $(X_j,\hor)$ as above ($j=1,\dots,d$), the group 
$
\prod\limits_{j=1}^d{\!_{\dps \hor}}\,\,\Aut(X_j,\hor)
$
acts on 
$
\prod\limits_{j=1}^d{\!_{\dps \hor}}\,\, X_j
$
by graph isometries via 
\begin{equation}\label{action}
gx = (g_1x_1)\cdots (g_dx_d)\,,\quad\text{where}\; g=g_1\cdots g_d
\AND x=x_1\cdots x_d\,.
\end{equation}

The basic example of a Busemann pair arises when the underlying graph is
a \emph{tree} $T$, that is, a connected graph without cycles, where 
$2 \le \deg(x) < \infty$ for every vertex~$x$. There are several
choices (one for each element in the boundary of the tree, see below)
to equip the edge set of $T$ with an orientation such that each vertex 
$x$ has a unique \emph{predecessor} $x^- $ and $\deg(x)-1$ successors $y \in T$ 
such that $y^-=x$. Then it is easily understood that in the induced half-order,
the \emph{ancestor relation} $\lle$, any two vertices $x, y \in T$
have a maximal common ancestor $x \cf y$. If $o \in T$ is a reference vertex (origin),
then we define $\hor(x) = d(x,o \cf x) - d(o, o \cf x)$, where $d(\cdot,\cdot)$
denotes the usual graph metric. Then $(T,\hor)$ is \emph{the} typical
example of a Busemann pair. 

In the present paper, we shall deal with \emph{homogeneous} trees  $T = \T_q$,
where each vertex has degree $q+1$ ($q \ge 2$). In this case, 
the horocyclic structure (i.e., the ancestor relation) is unique up to 
isomorphism. We write $\DL(q_1,\dots,q_d)$ for the horocyclic product
of the trees $T_1=\T_{q_1},\dots, T_d = \T_{q_d}$. The ``$\DL$'' stands for
{\sc Diestel and Leader}, who were the first~\cite{DiLe} 
to introduce the graph $\DL(2,3)$ in attempt to answer a question raised by 
{\sc Woess}~\cite{WoTop,SoWo}: ``is there a locally finite
vertex-transitive graph
which is not quasi-isometric with a Cayley graph of some finitely generated
group ?'' (Recall that a graph is called vertex-transitive if its isometry group
acts transitively on the vertex set.) A very recent announcement of
{\sc Eskin, Fisher and Whyte~\cite{EFW}} confirms that the graphs 
$\DL(q_1,q_2)$ ($q_1 \ne q_2$) are such examples.

The purpose of this paper is to present a picture of the many
interesting features of the graphs $\DL(q_1,\dots,q_d)$.

In \S \ref{iso}, we first recall in more detail the horocyclic structure of
the homogeneous tree $\T_q$ and the group $\Aff(\T_q) = \Aut(\T_q,\hor)$
of all its Busemann self-isometries. The latter group has been called the
\emph{affine group} of the tree by analogy with the affine  group over $\R$
acting on the hyperbolic upper half plane. We determine the full isometry
group $\Aut(\DL)$ of $\DL=\DL(q_1,\dots,q_d)$. We prove that
it is a finite extension of the group 
$
\Af = \prod\limits_{j=1}^d{\!_{\dps \hor}}\,\,\Aff(\T_{q_j})\,.
$
The latter acts transitively on $\DL$ and is \emph{amenable} as a 
locally compact, totally disconnected group with the topology of 
pointwise convergence. 

If the $q_j$ do not all coincide, we show that this group is also
\emph{non-unimodular} (i.e., the left Haar measure is not
right-invariant).  Consequently, by a theorem of {\sc Soardi and
  Woess}~\cite{SoWo}, the graph $\DL$ is \emph{non-amenable}, i.e., it
satisfies a strong isoperimetric inequality (the Cheeger inequality).
We also conclude that $\Aut(\DL)$ cannot have a co-compact lattice,
that is, there is no discrete (closed) subgroup that acts on $\DL$
with finitely many orbits. In particular, if the $q_j$ do not all
coincide, then $\DL$ is vertex-transitive, but is not the Cayley graph
of a finitely generated group.

In \S \ref{Cay}, we study $\DL_d(q) = \DL(q, \ldots, q)$, the
horocyclic product of $d$ copies of $\T_q$. We use an approach that is
reminiscent of the method for constructing lattices in Lie groups over
local fields, as outlined on the first page of the book by {\sc
  Margulis~\cite{Mar}}.  When $q=p_1\cdots p_r$ is the factorization
of $q$ as a product of prime powers, and $p_{\iota} \ge d-1$ for all
$\iota \in \{1,\dots,r\}$, the graph $\DL_d(q)$ is a Cayley graph of a
group of affine matrices over a ring of Laurent polynomials whose
coefficients come from a suitable finite ring. There is some degree of
freedom in the choice of the ring of coefficients. When $d=2$ or
$d=3$, we can take the ring $\Zq = \Z/(q\Z)$ of integers modulo $q$,
and for $d=2$ this is a way to describe the lamplighter group $\Zq \wr
\Z$, while for $d=3$ we obtain a finitely presented group into which
the lamplighter group embeds. This group has appeared in previous work
by {\sc Baumslag~\cite{Bau2}} and others.  In general, $\DL_d(q)$ is
\emph{quasi-isometric} with $\DL_d(q^s)$ for every $s \ge 1$. Thus,
$\DL_d(q)$ is always quasi-isometric with a Cayley graph of some
finitely generated group, while on the other hand, \cite{DiLe} and
\cite{EFW} suggest that the vertex-transitive graph
$\DL(q_1,\ldots,q_d)$ is not quasi-isometric with any Cayley graph
when the $q_j$ do not all coincide.

In \S \ref{sect:complex}, we consider $\DL(q_1,\dots,q_d)$ as a
$(d-1)$-dimensional cell complex and explore its homotopy type, which
is that of a union of countably many $(d-1)$-spheres glued together at
a single point.  This should be compared with a deep theorem of {\sc
  Bestvina and Brady~\cite{BeBr}}.  Thus, when $\DL$ is the Cayley
graph of a group, then this group is of type $F_{d-1}$, but not of
type $F_d$, and in particular it is finitely presented when $d \ge 3$.
We deduce that, for each $d$, the lamplighter group can be embedded in
a metabelian group of type $F_d$. In general, it is known~\cite{BoHa}
that every metabelian group embeds in a metabelian group of type
$F_3$, while embeddability in $F_d$ for larger $d$ is conjectured.

In \S \ref{sect:spectrum}, we turn our attention to a more
analytic-probabilistic object.  \emph{Simple random walk} on any
locally finite, connected graph $X$ is the Markov chain whose
transition matrix $P=\bigl(p(x,y)\bigr)_{x,y \in X}$ is given by
\begin{equation}\label{eq:SRW1}
p(x,y) = \begin{cases} 1/\deg(x)\,,&\text{if}\; y \sim x\,,\\
                       0\,,&\text{otherwise.}
\end{cases}
\end{equation}
$P$ acts on functions $f:X \to \R$ by 
\begin{equation}\label{eq:SRW2}
Pf(x) = \sum_y p(x,y)f(y)\,.
\end{equation}
In our case, $\deg(\cdot) = (d-1)(q_1+\dots+q_d)$ is constant, and we are 
interested in the \emph{spectrum} of $P$ on the space $\ell^2(\DL)$
of all square-summable complex functions on $\DL$. The spectral radius 
$\rho(P)$ is equal to~$1$ if and only if $q_1=\dots=q_d$. 
As a set, $\spec(P)$ is an interval contained in 
$[-\rho(P)/(d-1)\,,\,\rho(P)]$, and with the exception of a ``degenerate''
case, it coincides with the latter. In particular, for $\DL_d(q)$ the spectrum
of $P$ is the same as the spectrum of the projection of $P$ on the lattice
$\A_{d-1}$. The latter spectrum is absolutely continuous. On the other hand,
for arbitrary $q_1, \dots, q_d$, the spectrum of $P$ on $\DL(q_1, \dots, q_d)$
is \emph{pure point:} there is an orthonormal basis of $\ell^2(\DL)$ that 
consists of finitely supported eigenfunctions of $P$. This extends previous
results regarding the lamplighter group and the basic Diestel-Leader graphs
$\DL(q_1,q_2)$, see {\sc Grigorchuk and \.Zuk~\cite{GrZu2}}, 
{\sc Dicks and Schick~\cite{DiSc}} and {\sc Bartholdi and Woess~\cite{BaWo}}.
For the case $d=2$, the eigenvalues and eigenfunctions were computed
explicitly in those references. Here, we present explicit computations
for $d=3$ and $\DL_3(q)$, while the general case seems intractable
(except numerically).

Finally, in \S \ref{sect:Poisson}, we study the general class of random
walk on $\DL$ whose transition matrix is irreducible, invariant under the 
group $\Af$, and has finite first moment. Using results of
{\sc Cartwright, Kaimanovich and Woess~\cite{CaKaWo}} and 
{\sc Brofferio~\cite{Bro}} concerning random walks on $\Aff(\T_q)$, we show
that those random walks on $\DL$ converge almost surely to the geometric
boundary of $\DL$. The latter is the ideal boundary added to $\DL$
when considering the closure of $\DL$ in $\prod_i \wh T_i$, where $\wh T_i$
is the well-known \emph{end compactification} of $T_i$. We then use the
\emph{ray criterion} of {\sc Kaimanovich,} see \cite{KaWo}, to prove that
the active part of the boundary (i.e., the support of the limit distribution
of the random walk) is the ``largest possible'' model for distinguishing
limit points of the random walk: it is the \emph{Poisson boundary}.

This paper has become rather long and touches quite different aspects.
We have considered the possibility of splitting it in two parts, but
concluded that this would contradict its ``exploratory'' spirit.

\section{Isometry groups}\label{iso}

We start with a picture of the homogeneous tree $\T_2$ in horocyclic layers,
since it will be useful throughout the paper to keep this description
in mind. Note that the negative direction is ``upwards'' in the picture.

\vspace{-.4cm}

$$
\beginpicture 

\setcoordinatesystem units <.7mm,1.04mm>

\setplotarea x from -10 to 104, y from 4 to 84

\arrow <6pt> [.2,.67] from 2 2 to 80 80

\plot 32 32 62 2 /

 \plot 16 16 30 2 /

 \plot 48 16 34 2 /

 \plot 8 8 14 2 /

 \plot 24 8 18 2 /

 \plot 40 8 46 2 /

 \plot 56 8 50 2 /

 \plot 4 4 6 2 /

 \plot 12 4 10 2 /

 \plot 20 4 22 2 /

 \plot 28 4 26 2 /

 \plot 36 4 38 2 /

 \plot 44 4 42 2 /

 \plot 52 4 54 2 /

 \plot 60 4 58 2 /

 \plot 99 29 64 64 /

 \plot 66 2 96 32 /

 \plot 70 2 68 4 /

 \plot 74 2 76 4 /

 \plot 78 2 72 8 /

 \plot 82 2 88 8 /

 \plot 86 2 84 4 /

 \plot 90 2 92 4 /

 \plot 94 2 80 16 /


\setdots <3pt>

\putrule from -4.8 4 to 102 4
\putrule from -4.5 8 to 102 8
\putrule from -2 16 to 102 16
\putrule from -1.7 32 to 102 32
\putrule from -1.7 64 to 102 64
\setdashes <2pt>
\putrule from -1.7 -8 to 102 -8

\put {$\vdots$} at 32 -2
\put {$\vdots$} at 64 -2

\put {$\dots$} [l] at 103 6
\put {$\dots$} [l] at 103 48

\put {$H_{-3}$} [l] at -13 64
\put {$H_{-2}$} [l] at -13 32
\put {$H_{-1}$} [l] at -13 16
\put {$H_0$} [l] at -13 8
\put {$H_1$} [l] at -13 4
\put {$\bd^* \T$} [l] at -13 -8
\put {$\vdots$} at -10 -2
\put {$\vdots$} [B] at -10 70
\put {$\circ$} at 8 8
\put {$\omega$} at 82 82

\put {\scriptsize $0$} at 3.6 6.2
\put {\scriptsize $1$} at 12.2 6.2
\put {\scriptsize $0$} at 19.8 6.2
\put {\scriptsize $1$} at 28.4 6.2
\put {\scriptsize $0$} at 36   6.2
\put {\scriptsize $1$} at 44.2 6.2
\put {\scriptsize $0$} at 51.8 6.2
\put {\scriptsize $1$} at 60   6.2
\put {\scriptsize $0$} at 67.6 6.2
\put {\scriptsize $1$} at 76   6.2
\put {\scriptsize $0$} at 83.8 6.2
\put {\scriptsize $1$} at 92.2 6.2

\put {\scriptsize $0$} at 9 12
\put {\scriptsize $1$} at 22.5 12
\put {\scriptsize $0$} at 41 12
\put {\scriptsize $1$} at 54.5 12
\put {\scriptsize $0$} at 73 12
\put {\scriptsize $1$} at 86.5 12

\put {\scriptsize $0$} at 21 24
\put {\scriptsize $1$} at 43 24
\put {\scriptsize $0$} at 85 24

\put {\scriptsize $0$} at 45 48
\put {\scriptsize $1$} at 83 48

\put {\scriptsize $0$} at 72 75

\endpicture
$$

\vspace{.1cm}

\begin{center}
\centerline\emph{Figure 1}
\end{center}

\vspace{.1cm}

Along with that picture comes a more detailed description of the geometry
of $T=\T_q$. 

Any pair of vertices is connected by a unique 
\emph{geodesic path} $\geo{x\,y}$ whose length (number of edges) is the 
distance $d(x,y)$. A \emph{geodesic ray} is a one-sided infinite
geodesic path (isometric embedding of a half-line graph). Two rays are called
equivalent, if their symmetric difference (as sets of vertices) is finite.
An \emph{end} of $T$ is an equivalence class of rays. The boundary $\bd T$
of $T$ is the set of ends of $T$. For each $\xi \in \bd T$ and each $x \in T$
there is a unique geodesic ray $\geo{x\,\xi}$ that represents $\xi$ and starts
with $x$. We choose an origin (root) $o \in T$ and write $|x|=d(x,o)$.
If $v, w \in \wh T = T \cup \bd T$ then we can define their \emph{confluent} $c(v,w)$
as the last common element on $\geo{o\,w}$ and $\geo{o\,z}$, a vertex of $T$
unless $w=z \in \bd T$. With the ultrametric
$$
\theta(w,z) = \begin{cases} q^{-|c(w,z)|}\,,&\text{if}\;w\ne z\\
                            0\,,&\text{if}\;w = z\,,
	      \end{cases}	     
$$
$\wh T$ becomes a compact space. 

We now select an end $\om \in \bd T$ and write 
$\bd^* T = \bd T \setminus \{\om\}$. Given $\om$, we can define the
predecessor $x^-$ of $x \in T$ as the neighbour of $T$ that lies on $\geo{x\,\om}$.
Thus, the ancestor relation is 
\begin{equation}\label{ancestor}
x \lle y \Longleftrightarrow x \in \geo{y\,\om}\,,
\end{equation}
and for $x, y$ in general position, $x \cf y$ is the maximal common
ancestor, as explained in~\S 1. We write $\up(x,y) = d(x,x \cf y)$.
Then the horocycle index of $x$ (the Busemann function with respect 
to $\om$) is
$$
\hor(x) = \up(x,o) - \up(o,x)\,,
$$
and the $k$-th horocycle is $H_k = \{ x\in T : \hor(x)=k\}$. In particular,
\begin{equation}\label{updown}
d(x,y) = \up(x,y) + \up(y,x) \AND \hor(x) - \hor(y) = \up(x,y) - \up(y,x)\,.
\end{equation}

As in Figure~1, we can \emph{label} the edges of $T$ with the elements
of $\Zq$, such that the edges between a vertex and its $q$ successors carry
distinct labels, and such that on the geodesic from any vertex to $\om$,
only finitely many labels are non-zero. This labelling will be used several
times in the sequel.

For $T = \T_q$, its \emph{affine group} $\Aff(\T_q)$ is the stabilizer of
$\om$ in $\Aut(\T_q)$. It is an amenable and non-unimodular closed subgroup of 
$\Aut(\T_q)$ that acts transitively on $\T_q$, and all its elements
are Busemann isometries. We have $\hor(g) = \hor(go)$ for $g \in \Aff(\T_q)$.
See {\sc Cartwright, Kaimanovich and Woess}
\cite{CaKaWo} for more details about the structure of $\Aff(\T_q)$.

We shall need some basic facts about the modular function of an isometry group of 
a locally finite graph $X$ which is closed with respect to pointwise convergence. 
For more details, see {\sc Trofimov}~\cite{Tro} and
{\sc Woess}~\cite{WoTop}. If $\Ga \le \Aut(X)$ is such a group, and
$x \in X$, then $\Ga_x$ denotes the stabilizer of $x$ in $\Ga$, while
$\Ga x$ is the orbit of $x$ under $\Ga$. Since $\Ga$ is locally compact, it
carries a left Haar measure $dg$. The modular function $\Delta$ on $\Ga$ is the
unique multiplicative homomorphism $\Ga \to \R_+$ which satisfies
$$
\Delta(g_0) \int_{\Ga} f(gg_0)\,dg = \int_{\Ga} f(g)\,dg 
$$
for every $g_0 \in \Ga$ and every continuous, compactly supported function $f$ 
on $\Ga$. Inserting for $f$ the indicator function of $\Ga_x$ (which is
an open, compact subgroup of $\Ga$), one finds the formula
\begin{equation}\label{modular}
\Delta(g) = |\Ga_x(gx)|/|\Ga_{gx}x|
\end{equation}
for $g\in \Ga$ and for arbitrary $x \in X$, where $|\Ga_xy|$ is the
(finite) number of elements in the $\Ga_x$-orbit of $y$; see
e.g.~\cite{Tro,WoTop}. In particular, one has the following.

\begin{lem}\label{unimo}
If $\Ga$ acts transitively on $X$ then $\Ga$ is unimodular if and only
if $|\Ga_xy| = |\Ga_yx|$ for some ($\!\!\iff$ every)
$x \in X$ and all its neighbours $y$. 
\end{lem}

In the sequel, we fix integers $q_1, \dots, q_d \ge 2$, and write $o_j$ for the
origin of $T_j = \T_{q_j}$, while the symbol $o$ will be used for
the origin $o=o_1\cdots o_d$ of $\DL=\DL(q_1,\dots,q_d)$.
If $x, y \in \DL$, then we say that a neighbour $y$ of $x$ has
\emph{type}  $\ee_i-\ee_j$, if $y_i^- = x_i$ and $y_j=x_j^-$.
In this case, $x$ is a neighbour with type $\ee_j-\ee_i$ of $y$.
We write $N_{i,j}(x)$ for the set of neighbours with type $\ee_j-\ee_i$ of $x$.

\begin{pro}\label{pro:groupA}
The group
$$
\Af = \Af(q_1,\dots, q_d) = 
\prod\limits_{j=1}^d{\!_{\dps \hor}}\,\,\Aff(\T_{q_j})\,
$$
acts transitively on $\DL=\DL(q_1,\dots,q_d)$ by \eqref{action}. 
It is amenable. Furthermore,
$\Af$ is unimodular if and only if $q_1 = \dots = q_d$.
\end{pro}

\begin{proof} Let $x=x_1\cdots x_d$ be in $\DL$. Then there are 
$g_j \in \Aff(\T_{q_j})$ such that $g_jo_j=x_j$, $j=1,\dots,d$.
Setting $g=g_1\cdots g_d$ as in  \eqref{action}, we get 
$g \in \Af$, since $\sum_j \hor(g_j) = \sum_j \hor(x_j) = 0$.
Thus, $go=x$, and the action is transitive. Amenability of $\Af$
follows from the fact that it is a closed subgroup of the direct
product of the amenable groups $\Aff(\T_{q_j})$.

Regarding unimodularity, let $x$ be in $\DL$. By construction,
$\Af_x$ must map every neighbour $y$ of $x$ of type $\ee_j-\ee_i$ to a 
neighbour of $x$ of the same type, and every permutation of this
type can be achieved. Now $x$ has exactly $q_j$ neighbours of 
type $\ee_j-\ee_i$. Therefore, $|\Af_xy|=q_j$, and (by exchanging
$x \leftrightarrow y$ and $i \leftrightarrow j$) $|\Af_yx|=q_i$.
If we vary $i, j$ ($i \ne j$) and apply Lemma \ref{unimo}, then we
see that our group is unimodular if and only if all $q_j$ coincide.  
\end{proof}
 
Besides the elements of $\Af$, there may be further isometries of
$\DL$.  Let $\SF = \SF(q_1,\dots,q_d)$ be the group of all permutations $\si$
of $\{1,\dots,d\}$ such that $q_{\si(j)} = q_j$ for all $j$.
Then $\SF$ acts on $\DL$ by
\begin{equation}\label{permute}
\si x = x_{\si^{-1}(1)}\cdots x_{\si^{-1}(d)}\,, 
\end{equation}
that is, $\si$ permutes identical trees in the horocyclic 
product. Thus, $\SF$ also acts on $\Af$ by group automorphisms 
$(\si, g) \mapsto g^{\si} = \si g \si^{-1}$. 
If all $q_j$ are distinct, then $\SF(q_1,\dots,q_d)$ is of course trivial.
We shall prove the following.

\begin{thm}\label{fullgroup} The full isometry group of $\DL(q_1,\dots,q_d)$
is the semidirect product of $\SF$ with $\Af$ with respect to
the action $(\si, g) \mapsto g^{\si}$, 
$$
\Aut(\DL) = \SF \ltimes \Af\,.
$$
Thus, $\Aut(\DL)$ is amenable, and it is unimodular if and only if all $q_i$
coincide.
\end{thm}

For the proof, we need a description of the (graph-theoretical) 
\emph{link} $N(x)$ of a vertex 
$x \in \DL$, that is, the subgraph of $\DL$ spanned by the neighbours of
$x$. Under the graph homomorphism $\Hor: \DL \to \A_{d-1}$, where 
$\Hor(x) = \kk = \bigl(\hor(x_1),\dots,\hor(x_d)\bigr)$, the link $N(x)$ maps onto the
link $N(\kk)$ in the lattice $\A_{d-1}$. The latter link has 
$(d-1)$-cliques (complete graphs on $d-1$ vertices) as its building
blocks. Namely, for $i \in \{1, \dots, d\}$, write
$$
S_i^+(\kk) = \{\kk + \ee_i - \ee_j : j \ne i \} \AND
S_i^-(\kk) = \{\kk + \ee_j - \ee_i : j \ne i \}\,.
$$
Each of those spans a complete subgraph of $N(\kk)$. We have
$N(\kk) = \bigcup_{i=1}^d S_i^+(\kk)= \bigcup_{i=1}^d S_i^-(\kk)\bigr)\,$,
$$
S_i^+(\kk) \cap S_j^+(\kk) = S_i^-(\kk) \cap S_j^-(\kk) =  \emptyset\,,\AND
S_i^+(\kk) \cap S_j^-(\kk) = \{ \kk + \ee_i - \ee_j \} \quad (i \ne j)\,.
$$
We write $S_i^{\pm}(x)$ for the set of all points in $N(x)$ which are mapped
to $S_i^{\pm}(\kk)$ by $\Hor$. Note that $S_i^+(x) \cap S_j^-(x) = N_{i,j}(x)$.
The edges in $N(x)$ are as follows. 

(1) If $y, z \in S_i^+(x)$ then there are $j, k \ne i$ such that
$\Hor(y) = \kk+\ee_i - \ee_j$ and $\Hor(z) = \kk+\ee_i-\ee_k$. In this case,
there is an edge between $y$ and $z$ if and only if $j \ne k$ and $y_i = z_i$, in which case
$z$ is a neighbour of type $\ee_j - \ee_k$ of $y$ (i.e., $z_j^-=y_j$ and
$y_k^- = z_k$). Thus, the subgraph of $N(x)$ that is mapped onto an edge
in $S_i^+(\kk)$ is the graph $D(q_i,q_i)$ consisting of $q_i$
independent edges with their endpoints; see Figure 2a.

(2) If $y, z \in S_i^-(x)$ then there are $j, k \ne i$ such that
$\Hor(y) = \kk+\ee_j - \ee_i$ and $\Hor(z) = \kk+\ee_k-\ee_i$. In this situation,
there is an edge between $y$ and $z$ if and only if $j \ne k$; furthermore,
$z$ is a neighbour of type $\ee_k - \ee_j$ of $y$ (i.e., $z_k^-=y_k$ and
$y_j^- = z_j$). Thus, the subgraph of $N(x)$ that is mapped onto the edge
$[\kk+\ee_j - \ee_i,\kk+\ee_k-\ee_i]$ in $S_i^+(\kk)$ is the complete
bipartite graph $K(q_j,q_k)$; see Figure 2b.

$$
\beginpicture 

\setcoordinatesystem units <4mm,4mm>

\setplotarea x from -22 to 0, y from 0 to 5

\plot -22 3 -16 3 /

\plot -22 4  -16 4 /

\plot -22 5  -16 5 /

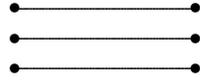
\put{Figure 2a: $D(3,3)$} [c] at -19 .9

\plot -6 3  0 3  -6 5  0 4  -6 3  0 5  -6 5 /

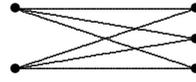
\put{Figure 2b: $K(2,3)$} [c] at -3 .9

\multiput {\scriptsize $\bullet$} at 
    -22 3  -16 3  -22 4  -16 4  -22 5  -16 5 
    -6 3  0 3  -6 5  0 4   0 5  /

\endpicture
$$

Figure 3a shows the link of a vertex of $\DL(2,2,3)$.
When $d \ge 3$, the link is connected. When $d=2$, it consists of
$q_1+q_2$ isolated points, and in this case, it will be more useful
to consider the $2$-link $N_2(x)$ spanned by all points at
distance $1$ and $2$ from $x$. Each of the $q_1$ neighbours $v$ 
of type $\ee_1-\ee_2$ of $x=x_1x_2$ 
is connected by an edge to each of the $q_2-1$ 
points $x_1y_2$, where $y_2 \ne x_2$ is a sibling of $x_2$ in
$T_2$, that is, $y_2^-=x_2^-$. In turn, there is an edge between 
each of those neighbours $v$ of $x$ and its $q_1$ neighbours 
of the same type $\ee_1-\ee_2$. Exchanging the role of $\ee_1$ and $\ee_2$,
one finds the other part of $N_2(x)$. See Figure 3b, where we have
also drawn the edges from $x$ to its neighbours in dotted lines.
$$
\beginpicture 

\setcoordinatesystem units <4mm,4mm>

\setplotarea x from -6 to 27.5, y from -6.5 to 6

\plot 4 0  2 3.464  -2 3.464  -4 0  -2 -3.464  2 -3.464  4 0 /

\plot 6 0  3 5.196  -3 5.196  -6 0  -3 -5.196  3 -5.196  6 0 /

\plot 2.5 4.33  -2.5 4.33 /

\plot 3 5.196  4 0  2.5 4.33  6 0  2 3.464 /

\plot -3 5.196  -4 0  -2.5 4.33  -6 0  -2 3.464 /

\plot -2 -3.464  3 -5.196 /

\plot 2 -3.464  -3 -5.196 /

\put {$S_1^+(x)$} [lt] at 4.5 -2.6
\put {$S_1^-(x)$} [rb] at -4.5 2.6
\put {$S_2^+(x)$} [rt] at -4.5 -2.6
\put {$S_2^-(x)$} [lb] at 4.5 2.6
\put {$S_3^+(x)$} [cb] at 0 5.3
\put {$S_3^-(x)$} [ct] at 0 -5.3

\put {Figure 3a: $N(x)$ in $\DL(2,2,3)$} [c] at 0 -8

\plot 16.5 5  17.5 3  18.5 5 /
\plot 19.5 5  20.5 3  21.5 5 /
\plot 16 0  17.5 3  19 0  20.5 3  16 0 /

\plot 19.5 -5  20.5 -3  25 0  26.5 -3  27.5 -5 / 
\plot 20.5 -5  20.5 -3  21.5 -5 /
\plot 25.5 -5  26.5 -3  26.5 -5 /
\plot 22.5 -5  23.5 -3  24.5 -5 /
\plot  23.5 -5   23.5 -3  25 0 / 

\setdots <3pt>
\plot 17.5 3  22 0  20.5 3 /
\plot 20.5 -3  22 0  23.5 -3 /
\plot 22 0  26.5 -3 /  
\setlinear

\put{$x$} [lb] at 22.2 0.2
\put {Figure 3b: $N_2(x)$ in $\DL(3,2)$} [c] at 20.5 -8

\multiput {\scriptsize $\bullet$} at 
    4 0  2 3.464  -2 3.464  -4 0  -2 -3.464  2 -3.464
    6 0  3 5.196  -3 5.196  -6 0  -3 -5.196  3 -5.196 
    2.5 4.33  -2.5 4.33 
    16 0  19 0  22 0  25 0
    17.5 3  20.5 3   
    16.5 5  18.5 5  19.5 5  21.5 5
    20.5 -3  23.5 -3  26.5 -3  
    19.5 -5  20.5 -5 21.5 -5  22.5 -5  23.5 -5  24.5 -5
              25.5 -5  26.5 -5  27.5 -5 /

\endpicture
$$

\begin{lem}\label{stab-permute} For every $x\in\DL$ and every $g$ in
  the stabilizer $\Aut(\DL)_x$ of $x$, there is
$\si \in \SF$ such that for all $i,j$ ($i \ne j$), we have 
$gN_{i,j}(x) = N_{\si^{-1}(i),\si^{-1}(j)}(x)$.
\end{lem}

\begin{proof} Our $g$ acts as a graph isometry on $N(x)$, and also on $N_2(x)$.
If $d=2$ (see Figure 3b) then $g$ must permute the two connected 
components of $N_2(x)$. This permutation must be trivial unless
$q_1=q_2$. Thus, the statement follows when $d=2$.

If $d \ge 3$, then by the above, each $S_i^-(x)$ is a complete $(d-1)$-partite 
subgraph of $N(x)$. That is, its vertex set is partitioned into
the sets $N_{j,i}(x)$ having cardinality $q_j$, $j \ne i$, such that
every pair of vertices in $N_{j,i}(x) \times N_{k,i}(x)$, $k \ne j$, is 
connected by an edge, while there are no edges between different vertices 
within each $N_{j,i}(x)$. Also, $S_i^-(x)$ is a \emph{maximal} complete
$(d-1)$-partite subgraph of $N(x)$ (i.e., it is not contained in any 
bigger complete $(d-1)$-partite subgraph).  
This property must be preserved by isometries of $N(x)$, and $gS_i^-(x)$
must again be a maximal complete $(d-1)$-partite subgraph of $N(x)$, whose
$d-1$ classes must have the same cardinalities $q_j$ as the sets $N_{i,j}(x)$,
$j \ne i$. Thus, there is a permutation $\si$ of $\{1\dots,d\}$ such that
$gS_i^-(x) = S_{\si^{-1}(i)}^-(x)$, and from
$|\{ j \ne i : q_j = q_i \} |= |\{ k \ne \si^{-1}(i) : q_k = q_i \}|$
we deduce $q_i = q_{\si^{-1}(i)}$ for each $i$. Therefore, $\si \in \SF$.

For each $j \ne i$, we must have $gN_{j,i}(x) = N_{k,\si^{-1}(i)}(x)$ 
 and $gN_{i,j}(x) = N_{l,\si^{-1}(j)}(x)$ for
some $k \ne \si^{-1}(i)$, $l \ne \si^{-1}(j)$. We still have to show that 
$k=\si^{-1}(j)$, and consequently $l=\si^{-1}(i)$.
Now note that with respect to its (``inner'') graph metric, $N(x)$ has 
diameter $3$, and that the only points at distance $3$ from 
all $y \in N_{i,j}(x)$ are precisely those in $N_{j,i}(x)$. Therefore
$$
N_{l,\si^{-1}(j)}(x) =  gN_{i,j}(x) = N_{\si^{-1}(i),k}(x)
$$
as sets, and thus $\bigl(l,\si^{-1}(j)\bigr) =\bigl(\si^{-1}(i),k\bigr)$.
This completes the proof of the lemma.
\end{proof}

We shall also need the following preparatory lemma.  

\begin{lem}\label{nij-permute} If $g \in \Aut(\DL)_x$ satisfies
$gN_{i,j}(x) = N_{i,j}(x)$ for all $i,j$ ($i \ne j$), then there 
is $h \in \Af_x$ such that $g\big|_{N(x)} \equiv h\big|_{N(x)}$.
\end{lem}

(As a matter of fact, it will turn out below that $g$ itself must belong
to $\Af$.)

\begin{proof} For each pair $(i,j)$, the isometry $g$ permutes the $q_i$
elements of $N_{i,j}(x)$. By the structure of $S_i^+(x)$, this permutation
must be independent of $j$ ($j \ne i$). That is, there is a permutation
$h_i$ among the successors of $x_i$ in $T_i$ such that $(gy)_i =
h_iy_i$ for all $y \in N_{i,j}(x)$ and all $j \ne i$. This permutation  
can be extended to an isometry of $T_i$, again denoted $h_i$, that fixes
$x_i$ and permutes the branches of $T_i$ ``below'' $x_i$. Setting 
$h = h_1 \cdots h_d$ according to \eqref{action}, we obtain the required
element of $\Af_x$.
\end{proof}

\begin{proof}[Proof of Theorem \ref{fullgroup}]
Consider an arbitrary $g \in \Aut(\DL)$. Set $x = g^{-1}o$, where
$o$ is the root of $\DL$. Then there is $g_x \in \Af$ such that $x = g_x^{-1}o$,
whence $gg_x^{-1}o = o$. By Lemma \ref{stab-permute}, there is $\si \in \SF$,
acting on $\DL$ by \eqref{permute}, such that 
$g' = \si^{-1}gg_x^{-1} \in \Aut(\DL)_o$ satisfies
$g' N_{i,j}(o) = N_{i,j}(o)$ for all $i, j$. We claim that $g'$
is \emph{type-preserving}, that is, for all 
$x \in \DL$,
\begin{equation}\label{type-preserve}
g' N_{i,j}(x) = N_{i,j}(g'x) \quad \forall\;i,j \in \{1,\dots,d\}
\; (i \ne j)\,.
\end{equation}
(Note that this does hold for every $g' \in \Af$.)
Since $\DL$ is connected and \eqref{type-preserve} is true for $x=o$, 
it is sufficient to show the following.
\begin{equation}\label{claim1}
\text{If \eqref{type-preserve} holds for $x \in \DL$ then it holds for 
every $y \in N(x)$.}
\end{equation}

So suppose that \eqref{type-preserve} holds for $x$. Let $h \in \Af_x$ be as in
Lemma \ref{nij-permute}, associated with $g'$, and set $g'' = h^{-1}g'$
Then $g''v = v$ for every $v \in \{ x \} \cup N(x)$. Let $y \in N(x)$,
so that $y \in N_{i,j}(x)$ for some $i, j$. Since $g''y=y$, Lemma 
\ref{stab-permute} implies that there is $\tau \in \SF$ such that
$g''N_{i,j}(y) = N_{\tau^{-1}(i),\tau^{-1}(j)}(y)$. 

Since $x \in N_{j,i}(y)$ and $g''x=x$, we find that $\tau(i) = i$
and $\tau(j)=j$. Also, if $k \ne i, j$, then $N_{k,i}(y)=N_{k,j}(x)$
is stabilized by $g''$. Thus, $\tau(k)=k$. We see that $\tau$ is the
identity, and as $h \in \Af$, 
$$
g' N_{i,j}(y) = hg'' N_{i,j}(y) = h N_{i,j}(y) = N_{i,j}(hy) = N_{i,j}(g'y)\,,
$$
as claimed. This proves \eqref{claim1} and consequently \eqref{type-preserve}.

We now use \eqref{type-preserve} to show that $g' \in \Af$.   
Our claim is the following.
\begin{equation}\label{claim2}
\text{If $x, y  \in \DL$ satisfy $x_i = y_i$ for some $i \in \{1,\dots,d\}$ 
then $(g'x)_i = (g'y)_i$.}
\end{equation}
Indeed, if this holds, then define $g_i\in\Aff(\T_{q_i})$ as follows.
Given $x_i \in T_i$, choose $x \in \DL$ with $i$-th coordinate $x_i$,
and set $g_ix_i = (g'x)_i$. This is independent of the specific choice of
$x$ by \eqref{claim2}. We therefore get $g'=g_1\cdots g_d \in \Af$.

For any $d$ and $i,j \in \{1, \dots, d\}$ with $i \ne j$, we can
define recursively $N_{i,j}^{(0)}(x) = \{ x \}$ and $N_{i,j}^{(k)}(x)
= \bigcup \{ N_{i,j}(y) : y \in N_{i,j}^{(k-1)}(x) \}$, where $x \in
\DL$. We observe that
$$v_j=w_j \quad\text{for all}\quad v,w \in N_{i,j}^{(k)}(x).$$

The proof of \eqref{claim2} is different in the cases $d=2$ and $d \ge
3$.  In both cases, we may assume without loss of generality that
$i=1$ in \eqref{claim2}. Suppose therefore that $x, y \in \DL$ satisfy
$x_1 = y_1$.  \smallskip

\underline{Case 1: $d = 2$.}\quad Consider $x_2, y_2$
and their common ancestor $u_2=x_2 \cf y_2$ in $T_2$. Since $x_1=y_1$, we
have $\hor(x_2)=\hor(y_2)$, whence 
$d(x_2,u_2) = d(y_2,u_2) = \hor(x_2) - \hor(u_2) =: k \ge 0$.
We can find $u_1 \in T_1$ such that $x_1 \lle u_1$ and  
$\hor(x_1) - \hor(u_1) = -k$, so that $u=u_1u_2 \in \DL$ and
$x,y \in N_{2,1}^{(k)}(u)$.  
Since $g'$  is type-preserving, we have 
$$
g'N_{2,1}^{(k)}(u) = N_{2,1}^{(k)}(g'u) \ni g'x, g'y\,.
$$
Using the above observation, we get $(g'x)_1 = (g'y)_1$, which proves
\eqref{claim2}.

\smallskip

\underline{Case 2: $d \ge 3$.}\quad
The subgraph of $\DL(q_1, \dots,q_d)$ spanned by 
$\{ v = v_1\cdots v_d \in \DL : v_1 = x_1\}$ is \emph{connected}; indeed, it is
isomorphic with $\DL(q_2, \dots, q_d)$.
Thus, there is a path from $x$ to $y$ in $\DL$ all whose vertices have
the same first coordinate $x_1$: if $v, w$ are successive vertices on this path
then $w \in N_{i,j}(v)$ where $i,j \ne 1$, and $v_1=w_1=x_1$. But then
\eqref{type-preserve} implies $g'w \in N_{i,j}(g'v)$, so that $g'v$ and $g'w$
differ only in the $i$-th and $j$-th coordinates. In particular,
$(g'w)_1 = (g'v)_1$, whence inductively $(g'y)_1=(g'x)_1$.

\smallskip

We conclude that $g  = \si g'$, where $\si \in \SF$ and $g' \in \Af$,
so that we have completed the description of $\Aut(\DL)$.

If $h \in \Af$ and $g=\si g' \in \Aut(\DL)$ with $\si \in \SF$
and $g' \in \Af$, then $ghg^{-1} = \si (g'hg'^{-1}) \si^{-1}$ is
type-preserving, so that \eqref{claim2} implies $ghg^{-1} \in \Af$. 
It is now obvious that the factor group $\Aut(\DL)/\Af$ is $\SF$.
\end{proof}

We briefly remind the reader of the concept of \emph{amenability}
of a locally compact group ($\equiv$ existence of a finitely additive,
left invariant probability measure on the group): see {\sc Paterson~\cite{Pat}}. Recall that a locally
finite graph $X$ is called amenable, if its \emph{isoperimetric constant}
is $0$, that is, the number
$$
\kappa = \inf \{ |\bd F|/\Vol(F) : F \subset X \;\text{finite}\,\}\,,
$$
where $|\cdot|$ is cardinality, $\bd F$ is the set of edges between
$F$ and $X \setminus F$, and $\Vol(F) = \sum_{x \in F} \deg(x)$.
By the well-known criterion of F\o lner, a finitely generated group is amenable
if and only if one ($\!\! \iff$ each) of its Cayley graphs is amenable in 
the latter sense.

\begin{cor}\label{isop}
The graph $\DL(q_1, \dots, q_d)$ is amenable if and only if all $q_i$
coincide.
\end{cor}

\begin{proof}
By a theorem of {\sc Soardi and Woess}~\cite{SoWo}, a vertex-transitive
graph is amenable if and only if its automorphism group (or equivalently,
any closed, transitive subgroup) is both amenable
and unimodular. Theorem \ref{fullgroup} yields the result.
\end{proof}

A \emph{co-compact lattice} in a locally compact group is a discrete,
countable subgroup such that the factor space is compact. 
The following is standard.

\begin{lem}\label{lattice}
If $X$ is a locally finite, vertex-transitive graph, then a group $\,\Ga$
is a co-compact lattice in $\Aut(X)$ if and only if $\,\Ga$ acts on $X$
with finitely many orbits and finite vertex stabilizers.
\end{lem}
 
\begin{proof} For safety's sake, we prove the ``only if''; the ``if'' is then
an easy exercise. If $\Ga$ is co-compact then there is a relatively compact
set $K \subset \Aut(X)$ such that $\Aut(X) = \Ga\, K$. If $o \in X$ is
arbitrarily chosen (a ``root''), then relative compactness of $K$ means that  
$F = Ko \subset X$ is finite. Transitivity of $\Aut(X)$ implies that 
$\Ga \,F = X$, whence $\Ga$ acts with finitely many orbits. 

Suppose that for some $x$ the stabilizer $\Ga_x$ contains a
sequence of distinct elements. Since the stabilizer of $x$ in $\Aut(X)$
is compact, that sequence must have an accumulation point. This contradicts
discreteness of $\Ga$.  
\end{proof}


%

A non-unimodular group cannot contain a co-compact lattice.
Recall that when $\Ga$ is a finitely generated group and $S=S^{-1} \subset \Ga$
is a finite set of generators, then the \emph{Cayley graph} $X(\Ga,S)$ of 
$\Ga$ with respect to $S$ has vertex set $\Ga$, and the edges are all
$[x,y]$, where $x, y \in \Ga$ and $x^{-1}y \in S$. The group $\Ga$ acts by
left multiplication on $X(\Ga,S)$, the action is transitive, and all vertex
stabilizers are trivial.
In view of Theorem \ref{fullgroup}, we now get the following.

\begin{cor}\label{nolattice-cor} If the $q_j$ do not all coincide, then
  $\Aut(\DL)$ does not possess any co-compact lattice. In particular,
  $\DL(q_1, \dots, q_d)$ cannot be the Cayley graph of a finitely
  generated group.
\end{cor}

\section{Cayley  graphs}\label{Cay}

We now study in more detail the case when $q_1= \dots = q_d = q$
all coincide. In this case, we write $\DL_d(q)$ for the horocyclic product of
$d$ trees with degree $q+1$.   
The fact that $\DL_2(q)=\DL(q,q)$ is a Cayley graph of the lamplighter group
$\Zq \wr \Z$, where $\Zq=\Z/(q\Z)$, has been explained (and exploited) 
in the papers by
{\sc Woess}~\cite{Wo} alone and with {\sc Bartholdi}~\cite{BaWo} and
{\sc Brofferio}~\cite{BroWo}. 
We now study the general case. 

Let $\Lq$ be a commutative ring of order $|\Lq|=q$ with multiplicative unit
$1$, and choose $\ell \in \Lq$.
Consider the ring $\Lq(\!(\tsf+\ell)\!)$ of all formal Laurent series
\begin{equation}\label{Lau-ser}
\fsf = \fsf(\tsf+\ell) = \sum_{k=-\infty}^{\infty} a_k\,(\tsf+\ell)^k\,,\quad 
a_k \in \Lq\,,\;\ \exists\ n_0 \in \Z: a_k = 0 \;\forall\ k < n_0\,.
\end{equation}
If $a_{n_0} \ne 0$ in this representation, then we write
$\vsf_{\ell}(\fsf) = n_0$ for the \emph{valuation}
of $\fsf$ with respect to $(\tsf+\ell)$, setting $\vsf_{\ell}(0) = \infty$.
The corresponding \emph{absolute value} is 
$$
|\fsf|_{\ell} = q^{-\vsf_{\ell}(\fsf)}\,.
$$
It defines a complete ultrametric on $\Lq(\!(\tsf+\ell)\!)$. Every closed ball 
in this  ultrametric has radius $q^n$ for some $n \in \Z$ and is also open. 

We shall write $\Lq(\!(\tsf+\ell)\!)^*$ for the multiplicative group of
\emph{invertible} elements in $\Lq(\!(\tsf+\ell)\!)$.

Following the explanation given in~\cite{CaKaWo}, we can identify
$\Lq(\!(\tsf+\ell)\!)$ with the ``lower'' boundary $\bd^* \T_q$, where
the vertex set of $\T_q$ is the set of all closed balls $x=
B_{\ell}(\fsf,q^n)$ with radius $q^n$, $n \in \Z$, and midpoint
$\fsf \in \Lq(\!(\tsf+\ell)\!)$. Relating this with the description of
\S 2, the predecessor of vertex $x$ as above is $x^- =
B_{\ell}(\fsf,q^{n+1})$.  The vertices on the $n$-th horocycle are, by
definition, precisely all balls with radius $q^{-n}$, which form a
partition of $\Lq(\!(\tsf+\ell)\!)$.

Conversely, starting with $\T_q$, for any $k \in \Z$, we can label the
edges between each vertex in $H_k$ and its $q$ successors with the
elements of $\Lq$ as in Figure 1, so that the sequence of labels on
the geodesic from $\om$ to $x\in H_n$
is $(a_k)_{k < n}$ with only finitely many $a_k \ne 0$.  Then $x=
B_{\ell}(\gsf,q^{-n})$, where $\gsf(\tsf) = \sum_{k < n}
a_k\,(\tsf+\ell)^k$.

Using this description, one sees that the
group $\Aff\bigl(\Lq(\!(\tsf+\ell)\!)\bigr)$ of all affine matrices
\begin{equation}\label{t-matrix}
g = \begin{pmatrix} \asf & \bsf \\ 0 & 1\end{pmatrix}\,,\qquad
\asf \in \Lq(\!(\tsf)\!)^*\,,\; \bsf \in \Lq(\!(\tsf)\!)
\end{equation}
embeds into 
$\Aff(\T_q)$ as a closed subgroup that acts transitively on the tree.
Namely, the image of a ball $x= B_{\ell}(\fsf,q^n)$ as 
above under such a matrix is the ball 
\begin{equation}\label{ball1}
gx = B\bigl(\bsf+\asf \cdot \fsf,q^{n-m}\bigr)\,,\quad\text{where}\quad
m=\vsf_{\ell}(\asf)\,.
\end{equation}

In the same way, we can also consider the ring $\Lq(\!(\tsf^{-1})\!)$ of all 
formal Laurent series $\fsf = \sum_k a_k\,\tsf^{-k}$ over $\Lq$ in the 
variable $\tsf^{-1}$ with the valuation $\vsf_{\infty}(\fsf)=n_0$
if $n_0$ is minimal such that $a_{n_0} \ne 0$. Again, we have the
identification with $\bd^* \T_q$, but in this specific case, 
we define the $n$-th horocycle to consist of all closed balls 
in $\Lq(\!(\tsf^{-1})\!)$ with radius $q^{-n-1}$. Again, the group
$\Aff\bigl(\Lq(\!(\tsf^{-1})\!)\bigr)$ embeds into 
$\Aff(\T_q)$ as a closed subgroup that acts transitively.

Now suppose that there are distinct elements 
$\ell_1, \ldots, \ell_{d-1} \in \Lq$ such that their pairwise differences
$\ell_i-\ell_j$ ($i \ne j$) are invertible. Consider the ring
\begin{equation}\label{eq:thering}
\Ring_d(\Lq) 
= \Lq[(\tsf+\ell_1)^{-1},\ldots,(\tsf+\ell_{d-1})^{-1},\tsf] 
\end{equation}
of all polynomials over $\Lq$ in the variables 
$(\tsf+\ell_1)^{-1},\ldots,(\tsf+\ell_{d-1})^{-1},\tsf\,$. 
They are of the form
\begin{equation}\label{eq:Lau-Pol}
\Psf=\Psf(\tsf+\ell_1,\ldots,\tsf+\ell_{d-1})
= \sum_{k_1, \dots, k_{d-1} \in \Z} a_{k_1, \dots, k_{d-1}}\,
(\tsf+\ell_1)^{k_1} \cdots (\tsf+\ell_{d-1})^{k_{d-1}}\,,
\end{equation}
where only finitely many of the coefficients 
$a_{k_1, \dots, k_{d-1}} \in \Lq\,$ are $\ne 0$.

Note that $\Ring_d(\Lq)$ depends not only on $d$ and $\Lq$, but also
on the choice of the $\ell_i\,$. The same is true for the group $\Ga_d(\Lq)$
in the next theorem.

\begin{thm}\label{thm:dtrees} 
If $\Lq$ contains distinct elements $\ell_1, \ldots, \ell_{d-1}$ such that 
their pairwise differences are invertible, then the Diestel-Leader graph
$\DL_d(q)$ is the Cayley graph
of the group $\Ga=\Ga_d(\Lq)$ of affine matrices
\begin{equation}\label{matrixd}
\begin{pmatrix} 
 (\tsf+\ell_1)^{k_1} \cdots (\tsf+\ell_{d-1})^{k_{d-1}} & \Psf\\ 
 0                                                      & 1
\end{pmatrix}\,,\qquad
k_1, \dots, k_{d-1} \in \Z\,,\; \Psf
                      \in \Ring_d(\Lq)
\end{equation}
with respect to the generators 
$$
\begin{aligned}
\begin{pmatrix} \tsf+\ell_i & b\\ 0 & 1\end{pmatrix}^{\textstyle \!\!\pm 1}\,,
\qquad\qquad
&b \in \Lq\,,\; i\in\{1,\ldots,d-1\} \AND\\
\begin{pmatrix} (\tsf+\ell_i)(\tsf+\ell_j)^{-1} & -b(\tsf+\ell_j)^{-1}\\ 
                 0 & 1\end{pmatrix}\,,\quad
                 &b \in \Lq\,,\; i,j\in \{1,\ldots,d-1\}\,,\;i \ne j\,. 
\end{aligned}
$$
\end{thm}

\begin{proof}
Each of $\Lq(\!(\tsf+\ell_i)\!)\;$ ($i = 1,\dots,d-1$) and  
$\Lq(\!(\tsf^{-1})\!)$ is the completion of $\Ring_d(\Lq)$ 
in its respective ultrametric. 
Indeed, for each $k \in \Z$ and distinct $i,j \in \{1,\dots,d-1\}$, 
we can write $(\tsf+\ell_j)^k$ as a Laurent series
in $(\tsf+\ell_i)$ as well as in $\tsf^{-1}$, namely
\begin{equation}\label{eq:changevar1}
(\tsf+\ell_j)^k = 
\sum_{n=0}^{\infty} \binom kn (\ell_j-\ell_i)^{k-n}(\tsf+\ell_i)^n =
\sum_{n=-k}^{\infty} \binom k{-n} \ell_j^{n+k} \tsf^{-n}\,.
\end{equation}
Also, for each \emph{non-negative} $k \in \Z$ and $j \in \{1,\dots,d-1\}$, 
we can write $\tsf^k$ as a polynomial in $\tsf+\ell_j\,$, 
\begin{equation}\label{eq:changevar2}
\tsf^k = \sum_{n=0}^{k} \binom kn(-\ell_j)^{k-n} (\tsf+\ell_j)^n\,.
\end{equation}
 
A straightforward computation shows that for all $k_1, \dots, k_{d-1} \in \Z$,
\begin{equation}\label{valuations}
\begin{aligned}
\vsf_{\ell_i}\bigl((\tsf+\ell_1)^{k_1} \cdots (\tsf+\ell_{d-1})^{k_{d-1}}\bigr) 
&= k_i\,,\quad i=1,\dots,d-1\,,
\AND\\
\vsf_{\infty}\bigl((\tsf+\ell_1)^{k_1} \cdots (\tsf+\ell_{d-1})^{k_{d-1}}\bigr)
& = -k_1-\dots -k_{d-1}\,,
\end{aligned}
\end{equation}
summing up to $0$. In particular, it follows from~\eqref{valuations} 
that the diagonal embedding
\begin{equation}\label{eq:diagembed}
\Ring_d(\Lq) \hookrightarrow \prod_{i=1}^{d-1}\Lq(\!(\tsf+\ell_i)\!) \times
\Lq(\!(\tsf^{-1})\!)\,,
\end{equation} 
where $\Psf \mapsto (\Psf_1, \dots, \Psf_d)$ with $\Psf_i = \Psf$ for all $i$,
is discrete in the latter direct product.  

Following the above explanations, we can identify the vertices
of the tree $\T_q$ with balls in the respective ultrametric. Thus,
we get $d$ copies of $\T_q$. The first $d-1$ correspond to 
$\Lq(\!(\tsf+\ell_i)\!)\equiv\bd^*\T_q\;$ ($i=1,\dots,d-1$), and recall
that in each of these cases, we define the $n$-th horocycle to consist
of all closed balls with radius $q^{-n}$ ($n\in \Z$).
The last copy of $\T_q$ corresponds to $\Lq(\!(\tsf^{-1})\!)\equiv\bd^*\T_q$, 
but also recall that this time, we define the $n$-th horocycle to consist 
of all closed balls with radius $q^{-n-1}$.

Each of the $d$ groups 
$\Aff\bigl(\Lq(\!(\tsf+\ell_i)\!)\bigr)$, $i=1,\dots,d-1$ and
$\Aff\bigl(\Lq(\!(\tsf^{-1})\!)\bigr)$ 
is a closed subgroup of $\Aff(\T_q)$ that acts transitively on the tree. 
By~\eqref{eq:diagembed}, the diagonal embedding of $\Ga_d(\Lq)$ into their 
direct product is discrete. By~\eqref{valuations}, combined with~\eqref{ball1}, 
the action of $\Ga_d(\Lq)$ on $(\T_q)^d$ is such that the subgraph $\DL_d(q)$ 
is preserved. The last step in showing that the latter is a Cayley graph of 
our group consists in exhibiting a correspondence between graph and group 
that is bijective and compatible with the group action:

\medskip

\noindent\emph{Claim.} For all $k_1, \dots, k_{d-1} \in \Z$, each  
$\Psf \in \Ring$ has a unique
decomposition into Laurent polynomials  
$\Psf = \Psf_{k_1, \dots, k_{d-1}}^{(1)} + \dots +
 \Psf_{k_1, \dots, k_{d-1}}^{(d)}$, where
\begin{equation}\label{decomp}
\begin{aligned}
\Psf_{k_1, \dots, k_{d-1}}^{(i)} &= 
\sum_{n < k_i} a_{n,i}\,(\tsf+\ell_i)^n\,,\quad i=1,\ldots,d-1\,, \AND\\
\Psf_{k_1, \dots, k_{d-1}}^{(d)} 
&= \sum_{n \ge -(k_1+\dots + k_{d-1})} a_{n,d}\,\tsf^n
\end{aligned}
\end{equation} 
with $a_{n,i} \in \Lq$ and $\ne 0$ only for finitely many $n$ ($i=1,\dots,d$).
\smallskip

\noindent\emph{Proof of the claim.} By multiplying with
$(\tsf+\ell_1)^{-k_1} \cdots (\tsf+\ell_{d-1})^{-k_{d-1}}$,
one sees that it is sufficient to prove this for $k_1= \dots k_d=0$. 
The existence of a decomposition~\eqref{decomp} is an easy 
exercise which we leave to the reader. Uniqueness requires a 
bit more care. To this end, we have to prove that 
\begin{equation}\label{eq:0decomp}
\Psf_{0, \dots, 0}^{(1)} + \dots + \Psf_{0, \dots, 0}^{(d)} = 0
\end{equation}
implies that $\Psf_{0, \dots, 0}^{(i)} = 0$ for all $i$.
Write $\Qsf$ for the left-hand side of~\eqref{eq:0decomp}. 
Suppose that for some $i \in \{1, \dots, d-1\}$, we have 
$\Psf_{0, \dots, 0}^{(i)} \ne 0$. Then 
$\vsf_{\ell_i}(\Psf_{0, \dots, 0}^{(i)})<\infty\,$, but by~\eqref{eq:changevar1}
and~\eqref{eq:changevar2}, 
$$
\vsf_{\ell_i}(\Psf_{0, \dots, 0}^{(i)}) = \vsf_{\ell_i}(\Qsf) 
= \vsf_{\ell_i}(0) = \infty\,,
$$
a contradiction. Thus, $\Psf_{0, \dots, 0}^{(i)} = 0$ for $i=1, \dots, d-1$,
and consequently also for $i=d$, concluding the proof of the claim. 

\medskip

Our choice for the origin is now $o=o_1 \cdots o_d$, where 
$o_i = B_{\ell_i}(0,1)\,$ for $i=1, \dots, d-1$ and
$o_d=B_{\infty}(0,q^{-1})$. 
Let 
$
g= \begin{pmatrix} 
 (\tsf+\ell_1)^{k_1} \cdots (\tsf+\ell_{d-1})^{k_{d-1}} & \Psf\\ 
 0 & 1
 \end{pmatrix} 
\in \Ga_d(q)\,.
$ 
In view of the claim we obtain a bijective correspondence 
$g \leftrightarrow x=x_1 \cdots x_d = go$, where 
$$
\begin{gathered}
x_i=g_io_i=B_{\ell_i}\bigl(\Psf_{k_1,\dots,k_{d-1}}^{(i)},q^{-k_i}\bigr)\,,
\quad i=1, \dots, d-1\,,\AND\\
x_d=g_do_d=B_{\infty}\bigl(\Psf_{k_1,\dots,k_{d-1}}^{(d)},
q^{k_1+ \dots + k_{d-1}-1}\bigr)\,.
\end{gathered}
$$
Here, when writing $go$, we intend the image of $g$ under diagonal 
embedding, while $g_io_i$ refers again to $g_i=g$, but this time
acting on the $i$-th tree of the horocyclic product: as isometries of
$T_i = \T_q$, the $g_i$ are distinct.
 
Compatibility with the group action is straightforward: 
$(gh)o = g(ho)$ for all $g, h \in~\Ga$. 

Finally, using \eqref{ball1}, one checks that for each $b \in \Lq$,
$\bigl(\begin{smallmatrix}\tsf+\ell_i & b\\0 & 1\end{smallmatrix}\bigr)$ maps
the origin $o$ to one of its neighbours of type $\ee_i - \ee_d$,
where $i \ne d$. Analogously, 
$\bigl(\begin{smallmatrix}(\tsf+\ell_i)(\tsf+\ell_j)^{-1} &
  -b(\tsf+\ell_j)^{-1}\\0 & 1\end{smallmatrix}\bigr)$ maps
$o$ to one of its neighbours of type $\ee_i - \ee_j$, where $i \ne j$.
\end{proof}

We now study how we can implement Theorem~\ref{thm:dtrees} for
specific values of $q$ and $d$.
\\[8pt]
{\bf A. $\DL_3(q)$ and a finitely presented extension of the
  lamplighter group.}
\\[2pt]
If $q \ge 2$ is arbitrary and $d=2$ (two trees) or $d=3$ (three
trees), then we can use for $\Lq$ the ring $\Zq$ of integers
modulo~$q$.  Indeed, we may choose $\ell_1=0$ and $\ell_2=1$.

Thus, for the lamplighter group, we start with the ring 
$\Ring_2(\Zq) = \Zq[\tsf^{-1},\tsf]$ of all Laurent polynomials 
$\Psf = \Psf(\tsf)= \sum_{k=K}^{L} a_k\,\tsf^k$ (with integers $K \le L $
and $a_k \in \Zq$).

\begin{cor}\label{cor:2trees}
The graph $\DL_2(q)=\DL(q,q)$ is a Cayley graph of the
Lamplighter group $\Zq \wr \Z$, represented as the group $\Ga_2(\Zq)$ of
affine matrices
$$
\begin{pmatrix} \tsf^k & \Psf
                  \\ 0 & 1\end{pmatrix}\,,\qquad
k\in \Z\,,\; \Psf
                      \in \Zq[\tsf^{-1},\tsf]
$$
with respect to the generators
$$
\begin{pmatrix} \tsf & b\\ 0 & 1\end{pmatrix}\,,
\quad b \in \Zq\,,
$$
and their inverses.
\end{cor}

For $\DL_3(q)$, we use $\Ring_3(\Zq) = \Zq[\tsf^{-1},(\tsf+1)^{-1},\tsf]$,
which consists of all
$$
\Psf=\Psf(\tsf,\tsf+1)= \sum_{k=K}^{L} \sum_{m=M}^{N}a_{k,m}\,\tsf^k (\tsf+1)^m\,,
$$ 
where $K, L, M, N \in \Z\,$, $\;K \le L\,$, $\;M \le N\,$,$\;a_{k,m} \in \Zq\,$.

\begin{cor}\label{cor:3trees} 
The Diestel-Leader graph $\DL_3(q) = \DL(q,q,q)$ is the Cayley graph
of the group $\Ga=\Ga_3(\Zq)$ of affine matrices
$$
\begin{pmatrix} \tsf^k(\tsf+1)^m & \Psf
                  \\ 0 & 1\end{pmatrix}\,,\qquad
k, m \in \Z\,,\; \Psf
                      \in \Zq[\tsf^{-1},(\tsf+1)^{-1},\tsf]
$$
with respect to the generators
$$
\begin{pmatrix} \tsf & b\\ 0 & 1\end{pmatrix}\,,\; 
\begin{pmatrix} \tsf+1 & b\\ 0 & 1\end{pmatrix}\,,\; 
\begin{pmatrix} \tsf^{-1}+1 &  -b\,\tsf^{-1}\\ 0 & 1\end{pmatrix}\,,
\quad b \in \Zq\,,
$$
and their inverses.
\end{cor}

The group $\Ga_3(\Zq)$ is finitely presented, see 
{\sc Baumslag~\cite{Bau1,Bau2}} 
and also below in \S \ref{sect:complex}. 
Thus, by Corollaries~\ref{cor:2trees}
and~\ref{cor:3trees} we recover a concrete way of embedding
the infinitely presented group $\Zq \wr \Z$ into the finitely presented,
metabelian group $\Ga_3(q)$.

\smallskip

How can we extend the above  to $\DL_d(q)$ with larger $d$~?
If $d \ge 4$, we cannot always find $d-1$ elements $\ell_i \in \Zq$ 
whose pairwise differences are invertible.
We note the following cases where we still can work with $\Zq$.
\\[4pt]
(i) When $q$ is \emph{odd} and $d=4$ then we can use $\ell_1=0$, $\ell_2=1$
and $\ell_3=-1 \equiv q-1$ (modulo~$q$). Thus, 
$\Ring_4(\Zq)= \Zq[\tsf^{-1},(\tsf+1)^{-1},(\tsf-1)^{-1},\tsf]$,
and $\DL_4(q)$ is a Cayley graph of the associated group of affine matrices.
\\[4pt]
(ii) When $q$ is \emph{prime} and $d \le q+1$ then we can use 
$\ell_1=0, \dots, \ell_{d-1}= d-2$. Again, $\DL_d(q)$ is a Cayley
graph of the group of affine matrices $\Ga_d(\Zq)$ associated
with $\Ring_d(q)=\Zq[\tsf^{-1},(\tsf+1)^{-1},\ldots, (\tsf+d-2)^{-1},\tsf]$.

\medskip

The last case (where $\Zq=\F_q$, the field of order $q$) points
towards a more general answer.
\\[8pt]
{\bf B.  $\DL_d(q)\,$ and quasi-isometries.}
\\[2pt]
For arbitrary $q \ge 2$, let $q=p_1 \cdots p_r$ be its factorization
as a product of prime powers. (\emph{Caution:} the $p_\iota$ are
prime powers, not necessarily primes themselves~!) 
Suppose first that $p_\iota \ge d-1$ for all $\iota
\in \{1, \ldots, r\}$, and let $\F_{p_\iota}$ be the field of order
$p_\iota$. Now construct the ring 
\begin{equation}\label{eq:ringL}
\F_q = \F_{p_1} \times \dots \times \F_{p_r}
\end{equation}
with coordinate-wise addition and multiplication. 
Let $\ell_{\iota,1}\,, \ldots, \ell_{\iota,d-1}$ be distinct elements of 
each $\F_{p_\iota}$, and set 
$\ell_j = (\ell_{1,j}\,, \ldots, \ell_{r,j}) \in \F_q$
for $j=1, \ldots, d-1$. Then all the differences $\ell_i - \ell_j\;$ ($i \ne j$)
are invertible elements of $\F_q$.
Therefore, we can apply Theorem~\ref{thm:dtrees} with 
the ring
$$
\Ring_d(\F_q) = \F_q[(\tsf+\ell_1)^{-1}, \ldots, (\tsf+\ell_{d-1})^{-1}, \tsf]
$$
of all Laurent polynomials over $\F_q$ in 
$\tsf+\ell_1\,,\ldots, \tsf+\ell_{d-1}\,$. 

\begin{cor}\label{cor:dtreesCayley} 
  Let $q=p_1 \cdots p_r$ be the factorization of $q$ in prime powers,
  and suppose $p_\iota \ge d-1$ for all $\iota \in \{1, \ldots, r\}$.
  Then $\DL_d(q)$ is a Cayley graph of the group $\Ga=\Ga_d(\F_q)$ of
  affine matrices $$
\begin{pmatrix} 
  (\tsf+\ell_1)^{k_1}\cdots(\tsf+\ell_{d-1})^{k_{d-1}} & \Psf \\ 
  0 & 1
\end{pmatrix}\,,\qquad
k_1,\ldots, k_{d-1} \in \Z\,,\; \Psf
                      \in \Ring_d(\F_q)
$$
\end{cor}

When $d \ge p_\iota$ for some $\iota \in \{1, \ldots, r\}$, then we cannot 
use the above construction. However, we can show that  
$\DL_d(q)$ is quasi-isometric with $\DL_d(q^s)$ for arbitrary $s \ge 1$, and
when $s$ is sufficiently large, the latter \emph{is} a Cayley graph
by Corollary~\ref{cor:dtreesCayley}. 

Recall that a quasi-isometry between two metric spaces $(X_1,d_1)$ and 
$(X_2,d_2)$ is a mapping $\phi: X_1 \to X_2$ such that there are constants
$A > 0$, $B \ge 0$ with
\begin{equation}\label{eq:quasiiso}
\begin{aligned}
A^{-1}d_2(\phi x_1,\phi y_1) -B &\le d_1(x_1,y_1) 
\le A \, d_2(\phi x_1,\phi y_1) + B \AND\\
d(x_2,\phi X_1) &\le B \quad \text{for all}\; x_1,y_1 \in X_1\,,\; x_2 \in X_2\,.
\end{aligned}
\end{equation}
For connected graphs, quasi-isometry refers to the discrete graph metric,
and for finitely generated groups, it refers to the word metric, i.e.,
the graph metric of one of its Cayley graphs: 
for a finitely generated group $\Ga$, any two of its 
Cayley graphs with respect to finite, symmetric generating 
sets are quasi-isometric (with $B=0$, i.e., bi-Lipschitz).
Quasi-isometry of groups has been an object of intensive study in geometric
group theory since its introduction by {\sc Gromov}~\cite{Gro}.
For the following, recall the group $\Af(q_1,\ldots,q_d)$ defined in 
Proposition~\ref{pro:groupA}.

\begin{pro}\label{pro:embed} For arbitrary $q_1, \ldots, q_d \ge 2$ and 
$s \ge 2$, the graphs $\DL(q_1^s, \dots, q_d^s)$ and
$\DL(q_1,\ldots,q_d)$ are quasi-isometric.

The quasi-isometry has the following property: if 
$\Ga \le \Af(q_1, \dots, q_d)$, then the subgroup
$$
\Ga_s = \{ g= g_1\cdots g_d \in \Ga : \hor(g_i) \in s\Z\,,\;i=1,\ldots,d  \} 
$$
embeds into $\Af(q_1^s, \dots, q_d^s)$, and if $\Ga$ acts transitively
and/or with trivial stabilizers on $\DL(q_1, \dots, q_d)$,
then it also acts transitively and/or with trivial stabilizers on 
$\DL(q_1^s,\ldots,q_d^s)$. 
\end{pro}

\begin{proof}
In the tree $\T_{q_i}$, consider the horocycles $H_{sn}$, $n \in \Z$.
Their union becomes the vertex set of $\T_{q^s}$, if we introduce
new edges as follows: if $n \in \Z$ and $x \in H_{sn}$, then draw an edge
between $x$ and $x^{-s} \in H_{s(n-1)}$, where $x^{-k}$ is defined
recursively by $x^{-1} = x^-$ (the predecessor of $x$), and 
$x^{-k-1} = (x^{-k})-$. See Figure~4, where $q=2$ and $s=3$.

$$
\beginpicture 

\setcoordinatesystem units <3mm,5mm>

\setplotarea x from -6 to 36, y from -2 to 6

\plot  0 0   1 2   2 0 /
\plot  4 0   5 2   6 0 /
\plot  8 0   9 2  10 0 /
\plot 12 0  13 2  14 0 /

\plot  1 2   3 4   5 2 /
\plot  9 2  11 4  13 2 /

\plot  3 4   7 6  11 4 /

\arrow <6pt> [.2,.67] from 17 3 to 20 3

\plot 22 0  29 6  24 0 /
\plot 26 0  29 6  28 0 /
\plot 30 0  29 6  32 0 /
\plot 34 0  29 6  36 0 /

\multiput {\scriptsize $\bullet$} at 
    0 0  2 0  4 0  6 0  8 0  10 0  12 0  14 0
    22 0  24 0  26 0  28 0  30 0  32 0  34 0  36 0  
    7 6  29 6 /

\setdots <3pt>

\putrule from -3.5 0 to 17.5 0
\putrule from -1.6 6 to 17.5 6
\putrule from 19.5 0 to 37 0
\putrule from 19.5 6 to 37 6

\put {$H_{sn}$} [l] at -6 0
\put {$H_{s(n-1)}$} [l] at -6 6

\put {Figure 4} [c] at 15 -2

\endpicture
$$
In this way, by inverting the arrow in Figure~4, we get a mapping 
$\phi_i: \T_{q_i^s} \to \T_{q_i}$ which is one-to-one from the horocycle 
$H_n(\T_{q_i^s})$ onto  $H_{sn}(\T_{q_i})$ for each $n \in \Z$. Furthermore, 
\begin{equation}\label{eq:phiproperties}
d(\phi_i x, \phi_i y) = s\cdot d(x,y) \AND \hor(\phi_i x) = s \cdot \hor(x)
\quad \text{for all}\; x, y \in \T_{q_i^s}\,.
\end{equation}
(Here, $d(\cdot,\cdot)$ refers to the respective graph metric in each of 
the two involved graphs.) Clearly, $\phi_i$ is a quasi-isometry, and 
by~\eqref{eq:phiproperties}, we can embed the group
$\{ g \in \Aff(\T_{q_i}) : \hor(g) \in s\Z \}$ into $\Aff(\T_{q_i^s})$
by $g \mapsto \phi_i^{-1} \circ g \circ \phi_i$.

We now define $\phi: \DL(q_1^s, \dots, q_d^s) \to \DL(q_1,\ldots,q_d)$
by $\phi = \phi_1\cdots \phi_d$, that is, the action of $\phi$ on the
$i$-th coordinate $x_i$ of $x=x_1\cdots x_d \in \DL(q_1^s, \dots, q_d^s)$
is given by  $x_i \mapsto \phi_ix_i$ ($i=1,\dots,d$).
It is now straightforward that $\phi$ is a quasi-isometry with the
asserted properties.
\end{proof}

\begin{cor}\label{cor:quasiiso}
For arbitrary $q$ and $d \ge 2$, the graph $\DL_d(q)$ is quasi-isometric
with $\DL_d(q^s)$, which is a Cayley graph of the group $\Ga_d(\F_{q^s})$,
when $s$ is sufficiently large.
\end{cor}

\begin{proof}
As above, let $q=p_1 \cdots p_r$ be the factorization of $q$ in
prime powers. Now choose $s$ sufficiently large such that
$p_\iota^s \ge d-1$ for all $\iota \in \{1, \ldots, r\}$. 
Then $\DL_d(q^s)$ is the Cayley graph of the group $\Ga_d(\F_{q^s})$
by Corollary~\ref{cor:dtreesCayley}. 
\end{proof} 

\begin{rmk}\label{rem:quasi} (a) In conclusion, we are not able to prove 
that \emph{every} Diestel-Leader graph $\DL_d(q)$ is itself a Cayley graph. 
The first open case is $\DL_4(2)$, the horocyclic product of $4$ trees 
with the same branching number $2$ ($\equiv$ degree $3$).
\\[4pt]
(b) Proposition~\ref{pro:embed} leads to the
following question:
is it true that two Diestel-Leader graphs are quasi-isometric if and 
only if there are $r, s \ge 1$ such that (up to permutation of their 
``coordinate'' trees) they are of the form $\DL(q_1^r,\dots, q_d^r)$ and  
$\DL(q_1^s,\dots, q_d^s)$. For the case $d=2$, compare with 
{\sc Wortman~\cite[\S 4]{Wort}} and the announcement of
{\sc Eskin, Fisher and Whyte~\cite{EFW}}.
\\[4pt]
(c) The following generalizes the conjecture of {\sc Diestel and
Leader~\cite{DiLe}} that $\DL(2,3)$ is not quasi-isometric to any Cayley
graph: is it true that $\DL(q_1,\cdots,q_d)$ for $d \ge 2$ is quasi-isometric
with a Cayley graph of some finitely generated group if and only if
$q_1 = \ldots = q_d$~? For $d=2$, see once more \cite{EFW}.
\end{rmk}
 
\section{The $\DL$ complex}\label{sect:complex}

We compute homotopical properties of $\DL$ in this section, and
combine this information with the results of the previous section to
derive conclusions on groups acting on $\DL$.

We defined $\DL$ as a graph in the introduction; for our purposes, it
is now better to view $\DL$ as a cell complex (in which cells are
represented by subsets of the vertex set), whose $1$-skeleton is the
graph defined by \eqref{eq:horoprod1} and \eqref{nbhd}.

Let therefore $T_1=\T_{q_1},\dots,T_d=\T_{q_d}$ be our homogeneous
trees with their respective Busemann functions. We now turn 
$\DL=\DL(q_1, \dots, q_d)$ into a cell complex, whose vertex set $\DL^0$
($0$-skeleton) is given by \eqref{eq:horoprod1}. Now
choose for each $i\in\{1,\dots,d\}$ a subset $E_i\subset X_i$ of
cardinality $1$ or $2$; if $E_i=\{x,y\}$, then $x\sim y$ has to be an edge of
$T_i$. To this choice there corresponds an $s$-dimensional \emph{cell}
$$
\gamma(E_1,\dots,E_d)=
\{x_1\dots x_d\in\DL^0:\,x_i\in E_i\;\text{for all}\;i\}\,,
$$
whenever $|\gamma(E_1,\dots,E_d)|\ge2$, with
$s=|E_1|+\dots+|E_d|-d-1$. The \emph{faces} of $\gamma(E_1,\dots,E_d)$ are
all $\gamma(E_1,\dots,E_{j-1},F_j,E_{j+1},\dots,E_d)$ for all possible 
choices of $j$ and $F_j\subset E_j$ with $|F_j|=1$ and $|E_j|=2$.

The $1$-skeleton of $\DL$ is given by the condition~\eqref{nbhd}, and 
the dimension of the $\DL$ complex is $d-1$. If $d=2$ then there are no cells of
dimension $\ge 2$, and we recover the original construction of $\DL$ as
a graph.

\begin{rmk}\label{twotriangles} An cell is always
of the form $\gamma(E_1,\dots,E_d)$, where each $E_i$ is either $\{x_i\}$ 
or $\{x_i,y_i\}$ with $x_i\sim y_i$ and $\hor(y_i)=\hor(x_i)+1$.

If $d\ge3$ then besides the $1$-dimensional cells (the edges), there are 
also two kinds of $2$-dimensional cells:
there must be precisely three
indices $j<k<l$ such that $|E_j|=|E_k|=|E_l|=2$. 
The requirement that $|\gamma(E_1,\dots,E_d)|\ge2$ is met when 
$\sum\hor(x_i) \in \{-1, -2\}$.
If now $\sum\hor(x_i)=-1$, then $\gamma(E_1,\dots,E_d)$ is a triangle,
with corners $x_1\cdots x_{j-1}y_jx_{j+1}\cdots x_d$, 
$x_1\cdots x_{k-1}y_kx_{k+1}\cdots x_d$ and
$x_1\cdots x_{l-1}y_lx_{l+1}\cdots x_d$.  If $\sum\hor(x_i)=-2$, then
$\gamma(E_1,\dots,E_d)$ is also a triangle but of a different kind:
its corners have the form $x_1\cdots y_j\cdots y_k\cdots x_d\,$, 
$x_1\cdots y_j\cdots y_l\cdots x_d$ and $x_1\cdots y_k\cdots y_l\cdots x_d$.

If $d\ge4$, then in addition, there are also three kinds of $3$-dimensional 
cells: there is now a set $I$ of four indices
with $|E_i|=2$ if $i\in I$, an we must have $\sum\hor(x_i) \in \{-1,-2,-3\}$. 
If $\sum\hor(x_i)=-1$, then $\gamma(E_1,\dots,E_d)$ is a tetrahedron, 
spanned by the
vertices $z_1\dots z_d$ with $z_j=y_j$ for precisely one $j\in I$, and
all other $z_i=x_i$. If $\sum\hor(x_i)=-2$, then
$\gamma(E_1,\dots,E_d)$ is an octahedron, spanned by the vertices
$z_1\dots z_d$ with $z_j=y_j$ for precisely two $j\in I$ (there are
$6=\binom 42$ choices). If $\sum\hor(x_i)=-3$, then
$\gamma(E_1,\dots,E_d)$ is again a tetrahedron, spanned by the
vertices $z_1\dots z_d$ with $z_j=y_j$ for precisely three $j\in I$.

In general, if $d \ge s+1$ there are $s$ kinds of $s$ dimensional cells,
according to the values of $\sum\hor(x_i) \in \{-1,\dots,-s\}$,
and $|E_i| = 2$ for all $i$ in a set $I \subset\{1, \dots, d \}$ of 
cardinality $s+1$.
\end{rmk}

\begin{dfn}\label{octahedron}
  Choose $R\in\N$ and points $b=b_1\cdots b_d$, $t=t_1\cdots t_d$,
  $t'=t_1'\cdots t_d' \in \prod_i T_i$ such that for all
  $i\in\{1,\dots,d\}$ we have (a)$\;\ \sum_i\hor(b_i)=-R\,,\;$ 
  (b)$\;\ b_i=t_i \cf t_i'\;$ (the maximal common ancestor), and 
  (c)$\;\ \hor(t_i)=\hor(t'_i)=\hor(b_i)+R\,$.

  The \emph{octahedron} $\Oct_{t,t',b}$ is the subcomplex of $\DL$
  spanned by the vertices $x_1\dots x_d$ such that $x_i$ lies on the
  path $\geo{t_i\,t'_i}$ for all $i$.

  The ``downward'' geodesics in $T_i$ of length $R$ starting at $b_i$
  can be ordered lexicographically by the sequence (``word'') of
  length $R$ of the labels (in $\{0,\dots,q_i-1\}$) along their edges.
  The octahedron $\Oct_{t,t',b}$ is \emph{basic} if for all $i$
  the geodesic $\geo{b_i\,t'_i}$ is the immediate successor of
  $\geo{b_i\,t_i}$, that is, all labels coincide except the first
  ones, which differ by~$1$. See Figure~5.
\end{dfn}
\vspace{-.3cm}
$$
\beginpicture 

\setcoordinatesystem units <5mm,6.5mm>

\setplotarea x from -6 to 6, y from -2 to 5


\plot  -5 0   0 5    5 0 /

\plot  10 0  11 1  10 2   9 3  10 4  11 5  
       12 4  11 3  12 2  13 1  12 0 /
 
\multiput {\scriptsize $\bullet$} at  
       -5 0  -4 1  -3 2  -2 3  -1 4  0 5  1 4  2 3  3 2  4 1  5 0  
       10 0  11 1  10 2  9 3  10 4  11 5  12 4  11 3  12 2  13 1  12 0 /

\multiput {$b_i$} [b] at 
               0 5.2  11 5.2 /
\multiput {$t_i$} [t] at 
               -5 -0.2   10 -0.2 /
\multiput {$t_i'$} [t] at 
                5 -0.2  12 -0.2 /

\multiput {\scriptsize $0$} [br] at 10.4 0.5  12.4 0.5  
        9.4 3.5  11.4 3.5  10.4 4.5 /

\multiput {\scriptsize $1$} [bl] at 11.6 4.5  
        9.6 2.5  11.6 2.5  10.6 1.5  12.6 1.5 /

\put {basic:} [br] at 8.5 4

\endpicture
$$
\vspace{-1.2cm}

\begin{center}
\centerline\emph{Figure 5}
\end{center}

\vspace{.2cm}

If $d=3$, octahedra are usual octahedra: six extremal vertices and
eight triangular faces with side lengths $R$. Each face is subdivided
into basic triangles with side lengths $1$, that is, $2$-cells; see
Figure~6. For arbitrary $d$, octahedra have
$2d$ vertices and $2^d$ top-dimensional faces.  $$
\beginpicture

\setcoordinatesystem units <5mm,4.5mm>

\setplotarea x from -21 to 6, y from -6.3 to 5

\plot  0 6  -6 0  2 -2  0 6  6 0  2 -2  0 -6 /
\plot  -6 0  0 -6  6 0 /

\setdashes <3pt>

\plot -6 0  -2 2  0 6 /
\plot  0 -6  -2 2  6 0 /

\multiput {\scriptsize $\bullet$} at  
       -6 0  0 6  6 0  0 -6  2 -2  -2 2 /

\put {$b_1b_2t_3$} [b] at 0 6.3
\put {$b_1b_2t_3'$} [t] at 0 -6.3
\put {$b_1t_2b_3$} [r] at -6.3 0
\put {$b_1t_2'b_3$\qquad\quad} [l] at 6.3 0

%

\put {$t_1b_2b_3$} [rt] at 1.6 -1.9
\put {$t_1'b_2b_3$} [lb] at -1.6 1.9

\put {\parbox{5cm}{\baselineskip 15pt A basic octahedron in $\DL(q_1,q_2,q_3)$.
      The face spanned by $b_1b_2t_3$, $b_1t_2b_3$ and 
      $t_1b_2b_3$ consists of all $x_1x_2x_3 \in \DL$
      with $x_i \in \geo{b_i\,t_i}$.}} [l] at -21 0
       
\endpicture
$$

\vspace{-.3cm}

\begin{center}
\centerline\emph{Figure 6}
\end{center}

\vspace{.3cm}

\newcommand\lfr{\mathfrak l}

\begin{lem}\label{lem:octahedron}
  Let $\lfr_i,\lfr'_i$ be distinct geodesics from $\omega_i$ to the
  ``lower'' boundary $\bd^*T_i$ in $T_i$, and set $\mathcal
  P=(\lfr_1\cup\lfr'_1)\times\cdots\times(\lfr_d\cup\lfr'_d)\cap\DL$.

  Let $b_i$ be the bifurcation point of $\lfr_i$ and $\lfr'_i$; set
  $R=-\sum_i\hor(b_i)$; and let $t_i,t'_i$ be the points respectively
  on $\mathfrak l_i,\mathfrak l'_i$ with
  $\hor(t_i)=\hor(t'_i)=\hor(b_i)+R$.

  If $R>0$, then $\mathcal P$ retracts to the octahedron
  $\Oct_{t,t',b}$. If $R\le0$, then $\mathcal P$ is contractible.

  Furthermore, every octahedron is homeomorphic to a sphere of
  dimension $d-1$.
\end{lem}

\begin{proof}
  Assume first $R>0$. Then clearly $\Oct=\Oct_{t,t',b}$ is a subset
  of $\mathcal P$.  Consider a point $x=x_1\dots x_d$ in $\mathcal
  P\setminus\Oct$. We will retract it to $\Oct$.  Consider (1)~in
  increasing order the coordinates $i\in\{1,\dots,d\}$ such that
  $\hor(x_i)<\hor(b_i)$ and, at the same time, (2)~in decreasing
  order, the coordinates $j$ such that $\hor(x_j)>\hor(t_i)$.  Move
  $x$, at unit speed, down on coordinate $i$ and up on coordinate $j$,
  until $x_i=c_i$ or $x_j\in\{t_j,t'_j\}$; when this happens, move to
  the next $i$ or $j$. If there are no more $j$'s available, keep
  moving up on the selected coordinate $j$.

  This process defines a retraction, i.e., a continuous map
  $\rho:\mathcal P\times[0\,,\,\infty)\to \mathcal P$ with
  $\rho(x,0)=x$ and $\lim_{t\to\infty}\rho(x,t)\in\Oct$.

  If $R\le0$, then this process also produces a retraction, but now
  towards a single point in $\mathcal P$.

  The last claim of the lemma is clear: an octahedron is topologically
  the $d$-fold join\footnote{recall that the join of two spaces $X,Y$
    is $X\star Y=X\times[0,1]\times
    Y/\{(x,0,y)\sim(x,0,y'),\,(x,1,y)\sim(x',1,y)\}$; it is classical
    that the join of a $(n-1)$-dimensional sphere and an
    $(m-1)$-dimensional sphere is a $(m+n-1)$-dimensional sphere.} of
  a pair of points, that is, a $(d-1)$-dimensional sphere.
\end{proof}
  
The following  result bears a strong similarity with Theorem~4.1 of
{\sc Bestvina and Brady~\cite{BeBr}}. 
However, there does not seem to be a natural discrete ambient group 
acting on $\prod T_i$ that would be required for their theorem to apply
directly.

\begin{thm}\label{homology}
The Diestel-Leader graph $\DL(q_1,\dots,q_d)$ has the homotopy type
of a wedge of countably many $(d-1)$-spheres, that is, the topological space
obtained by glueing together all those spheres at one single point.  
These spheres correspond bijectively to basic octahedra.
\end{thm}

\begin{proof}
Embed each tree $T_i$ in the upper half plane ${\mathbb H}$ in such a way that 
$\om_i \equiv \im\infty\,$, the ``lower'' boundary $\bd^*T_i$ lies in 
$\R$, and such that its edges are ordered lexicographically,
i.e., the edges labeled $0,1,\dots,q_i-1$ appear in left-to-right
order below any vertex. Let $\Omega_i$ be the set of all geodesics in
$T_i$ going from $\omega_i$ to $\R$.

The tree $T_i$, topologically, may be obtained from the disjoint union
of the geodesics in $\Omega_i$ by glueing the latter along 
specified half-geodesics from a vertex to $\om_i$.
Therefore, the product $T_1\times\dots\times T_d$ may be viewed as a
disjoint union of hyperplanes glued along ``hyper-octants'', and is 
naturally embedded in ${\mathbb H}^d$.

Similarly, $\DL$ is obtained by considering the disjoint union of the
planes $\lfr_1\times\dots\times\lfr_d\cap\DL$ for all $\lfr_i\in
\Omega_i$, and glueing them along subspaces.  It is naturally
embedded in a codimension-one contractible subspace $W$ of ${\mathbb
  H}^d$.

The homotopy type of $\DL$, therefore, is that of a wedge of
``pieces'', which are glued together at one single point.  Each
``piece'' is specified by two lexicographically consecutive geodesics
$\lfr_i,\lfr'_i$ in each tree $T_i$, and is of the form ${\mathcal
  P}=(\lfr_1\cup\lfr'_1)\times\cdots\times
(\lfr_d\cup\lfr'_d)\cap\DL$. It naturally is a subspace of $W$. 
By Lemma~\ref{lem:octahedron}, this piece $\mathcal P$ is either
contractible, or retracts to an octahedron, which is pure by the
choice of $\lfr_i,\lfr'_i$. Furthermore all basic octahedra appear in
this way.

Finally, it is obvious that there are countably many basic octahedra,
since there are countably many choices of $b_i$, $t_i$ and $t'_i$.
\end{proof}

Recall that a group $\Gamma$ is of \emph{type $F_d$} if it is the
fundamental group of an aspherical cell complex (i.e.\ a
$K(\Gamma,1)$) whose $d$-skeleton is compact. (``Aspherical'' means that
it is a topological space whose fundamental cover is contractible.) 
In particular, ``type $F_1$'' means ``finitely generated'', and 
``type $F_2$'' means ``finitely presented''.

\begin{cor}\label{fp} If $\DL_d(q)$ is a Cayley graph of a group $\Gamma$, 
then $\Gamma$ is of type $F_{d-1}$, but not of type $F_d\,$. In particular, 
if $d \ge 2$, it is finitely generated, and if $d\ge3$, it is 
finitely presented.
\end{cor}

\begin{proof}
The cell complex $\Ga\backslash\DL$ is  compact, and by
Theorem~\ref{homology} it approximates up to dimension $d-1$ a
classifying space $K(\Gamma,1)$. This proves that $\Gamma$ is of
type $F_{d-1}$. On the other hand, $K(\Gamma,1)$ has as many
$d$-cells as there are $\Gamma$-orbits on basic octahedra, and there
are infinitely many of them, so $\Gamma$ is not of type $F_d\,$.
\end{proof}

In particular, the group $\Ga_3(\Zq)$ is finitely presented, and
contains as a subgroup the lamplighter group
$\Ga_2(\Zq)\cong\Zq\wr\Z$. 

Explicit embeddings of the lamplighter
group were already known: {\sc Johnson~\cite{Jo}} embeds it first in
the group $G$ of permutations of $\Z$ that are translations outside of
a finite set (as the subgroup generated by the transposition $(1,2)$ and
the translation $n\mapsto n+2$); and then naturally embeds $G$ in the
group $\widetilde G$ of permutations of a three-branched star which
are translations outside of a finite set; he finally exhibits a
presentation of $\widetilde G$. (This group does \emph{not} act 
transitively and with trivial stabilizers on $\DL_3(q)$.)

{\sc Baumslag~\cite{Bau1,Bau2}} proves that every finitely
generated metabelian group can be embedded in a finitely presented
metabelian group. His construction is quite explicit: the group into which
$\Zq \wr \Z$ embeds has presentation
\begin{equation}\label{eq:Bpresentation}
\Gamma=\langle a,s,t \mid \,a^q,[s,t],[a,a^t],a^{-t}aa^s\rangle\,.
\end{equation}
In fact, this presentation can be interpreted in terms of ``lamplighters''
quite naturally as follows. 

Consider the abelian group $\mathfrak{V}$ of all finitely supported
functions $\Z^2 \to \Zq$ (the \emph{configurations} of ``lamps'').
Thus, $\mathfrak{V}$ is generated by all point masses $\delta_{x,y}$
at position $(x,y) \in \Z^2$ (with $\delta_{x,y}(x,y)=1\in\Zq$ and $0$
elsewhere).  Let $\mathfrak{W}$ be the subgroup generated by all
elements of the form $\delta_{x,y}+\delta_{x+1,y}-\delta_{x,y+1}$,
where $(x,y) \in \Z^2$.  The group $\Z^2=\langle s,t\rangle$ acts by
translations on $\mathfrak{V}$ preserving $\mathfrak{W}\,$, and one
has $\Gamma \cong \langle s,t \rangle \ltimes
(\mathfrak{V}/\mathfrak{W})\,$, by identifying the element $a$ of the
presentation \eqref{eq:Bpresentation} with $\delta_{0,0}$, and more
generally $a^{s^xt^y}$ with $\delta_{x,y}$.

Using Corollary~\ref{cor:3trees}, it is easy to identify $\Gamma$ with
$\DL_3(q)$, under the correspondence
\[a\leftrightarrow\begin{pmatrix} 1 & 1\\ 0 & 1\end{pmatrix},\quad
s\leftrightarrow\begin{pmatrix} \tsf^{-1} & 0\\ 0 & 1\end{pmatrix},\quad
t\leftrightarrow\begin{pmatrix} (1+\tsf)^{-1} & 0\\ 0 & 1\end{pmatrix}.\]

We also remark that in our context it may be more adequate to replace $\Z^2$
with $\A_2$ via the correspondence $(x,y) \mapsto (x,y,-x-y)$.

More generally, a presentation of $\Ga_d(q)$ may be read from $\DL$,
as follows.

\begin{thm}\label{thm:presentation}
Under the hypotheses and with the notation of Corollary \ref{cor:dtreesCayley}, if $d\ge 3$, 
the group $\Ga_d(q)$ admits the presentation
$$
\begin{aligned}
\Ga_d(q)=\bigl\langle &g_{i,j,\lambda} \,,\; 
  1\le i,j\le d\; (i \ne j)\,,\;\lambda\in\Zq \mid 
   g_{i,j,\lambda}\, g_{j,i,-\lambda}\,,\\
  &g_{j,i,\lambda\,}g_{k,j,\mu}\,g_{i,k,\nu}\;\text{whenever}\;\lambda+\mu+\nu=0,\\
  &g_{i,j,\lambda}\,g_{j,k,\mu}\,g_{k,i,\nu}\;\text{whenever}\;
   \lambda+\mu+\nu=0\;\text{and}\;\lambda\,\ell_k+\mu\,\ell_i+\nu\,\ell_j=0\;
   \text{and}\;d\neq i,j,k,\\
  &g_{i,j,\lambda}\,g_{j,d,\mu}\,g_{d,i,\nu}\;\text{whenever}\;\mu+\nu=0\;
   \text{and}\;\lambda+\mu\,\ell_i+\nu\,\ell_j=0\bigr\rangle.
\end{aligned}
$$
\end{thm}
\begin{proof}
By Theorem~\ref{homology}, the $2$-skeleton of $\DL$ is simply
connected if $d\ge 3$. Therefore a presentation of $\Ga_d(q)$ may be
read from this $2$-skeleton. The generators are the
$\Ga_d(q)$-orbits of edges, and the relations are the
$\Ga_d(q)$-orbits of $2$-cells. These orbits can be identified with
the neighbourhood of a given vertex. There are $d(d-1)q$ generators,
since the latter number is the vertex degree in $\DL_d(q)$.

The group element $g_{i,j,\lambda}$ maps $o$ to its neighbour $x$ of type
$\ee_i - \ee_j$ with label $\lambda$, that is, $x_j = o_j^-$, $x_i^- = o_i$,
and the label on the edge $[o_i,x_i]$ of $T_i$ is $\lambda$.
Thus, the inverse of $g_{i,j,\lambda}$ is $g_{j,i,-\lambda}$.

There are $d(d-1)(d-2)q^2$
relations of length $3$ corresponding to the first kind of
triangles indicated in Remark~\ref{twotriangles}; these are given by
the relations in the third line. There are also $d(d-1)(d-2)q$
relations of length $3$ corresponding to the second kind of
triangles indicated in Remark~\ref{twotriangles}; these are given by
the relations in the last two lines.
\end{proof}

Note that the generators $g_{i,d,\lambda}$ and $g_{i,j,\lambda}$ seem
to have different roles in the above presentation.  This asymmetry can
be masked as follows: define polynomials $\Lambda_i=\tsf+\ell_i$ for
$i<d$, and $\Lambda_d=1$. Then the presentation can be written as
$$
\begin{aligned}
\Ga_d(q)=\bigl\langle &g_{i,j,\lambda} \,,\; 
  1\le i,j\le d\; (i \ne j)\,,\;\lambda\in\Zq \mid 
   g_{i,j,\lambda}\, g_{j,i,-\lambda}\,,\\
  &g_{j,i,\lambda\,}g_{k,j,\mu}\,g_{i,k,\nu}\;\text{whenever}\;\lambda+\mu+\nu=0,\\
  &g_{i,j,\lambda}\,g_{j,k,\mu}\,g_{k,i,\nu}\;\text{whenever}\;
\lambda\Lambda_k+\mu\Lambda_i+\nu\Lambda_j=0\bigr\rangle.
\end{aligned}
$$

Moreover, combining with Corollary~\ref{fp}, we obtain:
\begin{cor}\label{cor:lamp-embed}
  The lamplighter group $\Zq\wr\Z$ can, for all $d$, be embedded
  in a metabelian group of type $F_d$.
\end{cor}

\begin{proof}
  We will use some algebraic number theory in the proof.

  Factor $q=p_1\cdots p_r$ as a product of prime powers. It suffices
  to embed $\Zyk_{p_\iota}\wr\Z$ in a metabelian group of type $F_d$.
  Indeed, $\Zq \cong \Zyk_{p_1} \times \ldots \times \Zyk_{p_r}$, so that by 
  identifying $\Z$ with the principal diagonal of $\Z^r$ via
  $k \leftrightarrow (k,\dots,k)$, we can embed $\Zq \wr \Z$ into the direct 
  product of the $\Zyk_{p_\iota}\wr\Z$. But type $F_d$ is inherited by 
  products.
  
  Write therefore $p_\iota=p^n$ as a prime power, and let $s$ be such
  that $p^s\ge d$.  The field $\F_{p^s}$ contains elements
  $\ell'_1,\dots,\ell'_d$ whose pairwise differences are invertible.

  There exists a unique non-ramified extension $\mathbb K$ of $\Q_p$, with
  ring of integers $\Int $ and maximal ideal 
  $\mathcal M=p\Int $, such that $\Int /\mathcal M\cong \F_{p^s}$.
  Furthermore, $\mathbb K$ is a $\Q_p$-vector space of dimension $s$, and
  $\Int $ is a free $\Z_p$-module of rank $s$. For all $i$, let
  $\ell''_i$ be an arbitrary preimage of $\ell'_i$ in $\Int $ under the
  natural projection $\Int  \to \F_{p^s}$. Then
  the pairwise differences $\ell''_i-\ell''_j$ lie in 
  $\Int \setminus\mathcal M$, and are therefore invertible in~$\Int $.

  Now set $\Lb=\Int /\mathcal M^n$; since
  $\Z_p/p^n\Z_p\cong\Zyk_{p^n}$, this is a free $\Zyk_{p^n}$-module.
  Let $\ell_i$ be the image of $\ell''_i$ in $\Lb$. The pairwise
  differences of the different $\ell_i$ are again invertible in $\Lb$.

  In the group $\Ga_{d+1}(\Lb)$ --- which has type $F_d$ by Corollary~\ref{fp}
  --- consider the elements
  $x=(\begin{smallmatrix}1&1\\0&1\end{smallmatrix})$ and
  $y=(\begin{smallmatrix}\tsf+\ell_1&0\\0&1\end{smallmatrix})$.
  Then $x$ has order $p^n$, and the subgroup $\langle x,y\rangle$
  of $\Ga_{d+1}(\Lb)$ is isomorphic to $\Zyk_{p^n}\wr\Z$.
\end{proof}

Note that it is not possible to embed $\Zq\wr\Z$ in a metabelian group
of type $F$ (i.e., one with a compact $K(\Gamma,1)$). Indeed recall
the following from {\sc Bieri and Groves~\cite{BiGr}}: a group $G$ is
of type $F\!P_n$, for $n\in\N\cup\{\infty\}$, if the trivial
$G$-module $\Z$ admits a $\Z G$-projective resolution of bounded rank
in each dimension up to $n$; this condition is weaker than type $F_n$.
A group $G$ is of type $F\!P$ if it is of type $F\!P_\infty$ and has
finite cohomological dimension; this condition is weaker than type
$F$, but it also implies that $G$ is torsion-free.  In \cite{BiGr} it is
proven that a metabelian group of type $F\!P_\infty$ is
virtually of type $F\!P$, and in particular is virtually torsion-free.

\begin{rmk}\label{remark:automat}
Assume that $\Gamma$ is one of the groups with Cayley graph $\DL_d(q)$, 
constructed in \S\ref{Cay}. Then $\Gamma$ is in most cases an
``automata group'', as we shall explain below. But it is never an
``automatic group'' in the sense of
{\sc Epstein et al.}~\cite{Eps+}. Indeed, by
Corollary~\ref{fp}, $\Gamma$ is of type $F_{d-1}$ but not of type $F_d$; in
particular it is not of type $F\!P_\infty$. On the other hand, automatic
groups are always of type $F\!P_\infty$,
by~\cite[Theorem~10.2.6]{Eps+}.

We briefly recall the definition of automata groups. Fix an alphabet
$\Sigma$. Automata groups are permutation groups of the set of
infinite words $\Sigma^{\infty}$. A family of \emph{automatic
  transformations} is given by a machine with bounded memory
computing the permutation in real time. It can be modeled by a finite
set of states $A$ and a function $\Phi:A\times\Sigma\to\Sigma\times A$
such that for all $a\in A$ the composition
$((\tau,b)\mapsto\tau)\circ\Phi\circ(\sigma\mapsto(a,\sigma))$ is a
permutation of $\Sigma$. The associated transformations $T_a$, for $a
\in A$, act on infinite strings $\sigma=\sigma_1\sigma_2\dots \in
\Sigma^{\infty}$ by
$$
T_a(\sigma) = \tau_1 T_b(\sigma_2\sigma_3\dots),
\text{ if }\Phi(a,\sigma_1)=(\tau_1,b).
$$
The product and inverse of automatic transformations are easily seen
to be automatic.

Assume now that the  group $\Gamma = \Ga_d(\Lq)$ is constructed as in 
Theorem~\ref{thm:dtrees} 
by use of elements $\ell_1,\dots,\ell_{d-1} \in \Lq$ whose pairwise
differences are invertible. \emph{Assume furthermore that
all $\ell_i\in\Lq$ are invertible.} Identify $\Sigma$ with $\Lq$, and
$\Sigma^{\infty}$ with $\Lq[[\tsf]]$, the ring of formal power
series in $\tsf$ over $\Lq$.

For $j\in\{1,\dots,d-1\}$, consider the finite set of states $A^{(j)}
\cong \Lq$ with associated affine transformations $T_a^{(j)}$ of
$\Sigma^{\infty} \equiv \Lq[[\tsf]]$ of the form
\begin{equation}\label{eq:transformation}
T_a^{(j)}\bigl(\fsf(\tsf)\bigr) = a+(\tsf+\ell_j)\fsf(\tsf)\,,\quad 
\text{where} \quad a\in\Lq\,.
\end{equation}
Note here that invertibility of $\ell_j$ implies that also the inverse
transformation of $T_a^{(j)}$ preserves $\Lq[[\tsf]]$.  If we write
$\fsf(\tsf) = b_0+\tsf\,\gsf(\tsf)$ with constant term $b_0 \in \Lq$
and $\gsf(\tsf) \in \Lq[[\tsf]]$, then $$
T_a^{(j)}\bigl(\fsf(\tsf)\bigr)=(\ell_j b_0+a) + \tsf \cdot
T_{b_0}^{(j)}\bigl(\gsf(\tsf)\bigr)\,, $$
and therefore each
$T_a^{(j)}$ is an automatic transformation, with $$
\Phi_j(a,b)=(\ell_j b + a,b) \,.  $$
Since $\Gamma$ is generated by
$\{ T_a^{(j)} : j=1, \dots, d-1\,,\; a \in \Lq \}\,$, all its elements
are automatic transformations.  We note again that in the construction
above we need $\ell_j$ to be invertible; the inverse of $T_a^{(j)}$ is
the transformation $T_a^{(j)}$ defined on the set of states
$A^{(j)}\cong\Lq$ by the map $\overline\Phi(a,b)=
\big((b-a)\ell_j^{-1},(b-a)\ell_j^{-1}\big)$.

For more general information on automata groups, see the 
survey by {\sc Bartholdi, Grigorchuk and Nekrashevych~\cite{BaGrNe}}. 
\end{rmk}

We finally remark here that $\DL(q_1,\ldots,q_d)$ has rational growth
function $ \sum_{x \in \DL} \tsf^{d(x,o)}\,, $ where $d(\cdot,\cdot)$
is the graph metric. We intend to come back to this point on another
occasion.

\section{The spectrum of simple random walk}\label{sect:spectrum}

Recall the definition \eqref{eq:SRW1}, \eqref{eq:SRW2} of the SRW 
operator $P$. It acts on the space $\ell^2(\DL)$ of all
square-summable functions $f:\DL \to \C$ with the standard inner
product $\langle f,g \rangle = \sum_x f(x)\overline{g(x)}$.
Since $P$ is stochastic and self-adjoint, we have $\|P\| = \rho(P) \le 1$,
where $\rho(P)$ is its spectral radius. We also set $\rho'(P) = \min \spec(P)$,
and write $D=(d-1)\sum_i q_i$ for the (constant) vertex degree of $\DL$.
\\[6pt]
{\bf A. Polyhedra and horizontal functions.}
\\[2pt]
We now generalize the method of~\cite{BaWo}, where the spectrum
of SRW on the horocyclic product of two trees was considered. 

A function $f:\DL \rightarrow \R$ is called \emph{horizontal\/} if 
it is finitely supported and 
\begin{equation}\label{eq:horizontal}
\sum_{y_j \in \T_j\,: \,\hor(y_j)=\hor(x_j)}  
f(x_1\cdots x_{j-1}y_jx_{j+1}\cdots x_d) = 0 
\end{equation}
for every $x=x_1\cdots x_d\in \DL$ and $j \in \{1,\dots,d\}\,$.

\begin{lem}\label{lem:dense}
The subspace of $\ell^2(\DL)$ spanned by the horizontal functions is dense.
\end{lem}

\begin{proof}
Since the point masses $\delta _x$ ($x \in \DL$) generate a dense 
subspace of $\ell^2(\DL)$, and since $\DL$ is vertex-transitive, it is
sufficient to show that $\delta _{o}$ (with $o=o_{1}\cdots o_{d}$) can be 
approximated in the $\ell^2$-norm by horizontal functions. 

For each $j$, pick a vertex $b_{j}=b_{j}^{n}$ of $T_{j}$ on the
horocycle $H_{-n}^{j}$ of $T_j$ that is not an ancestor of $o_j$.
Now let $f_{j}=f_{j}^{(n)}$ be the function on $T_j$ defined by 
\begin{equation*}
f_j (x_j) = \begin{cases} 1\,, &\text{if}\; x_j=o_j\,,\\
                -q_{j}^{-n}\,, &\text{if}\; b_j \lle x_j \in H_0^j\,, \\ 
                             0 &\text{in all other cases.}
            \end{cases}
\end{equation*}
Let $f=f^{(n)}$ be defined by 
$f(x_1\cdots x_d) = f_1(x_1) \cdots f_d(x_d)\,$. 
This function is horizontal, and  
\begin{equation*}
\left\| f^{(n)}-\delta _{o}\right\| ^{2}
=\prod_{j=1}^d (q_j^{-n}+1) -1 \rightarrow 0\,, \quad \text{as}\quad
n\to \infty\,. \qedhere
\end{equation*}
\end{proof}

\begin{dfn}\label{def:poly}
Let $a_{1},\ldots ,a_{d}$ be vertices of $T_{1},\ldots ,T_{d}$, respectively,
such that $\he=-\sum_{j=1}^{d}\hor( a_j) \geq 2$. The 
\emph{polyhedron\/} $\Pol=\Pol(a_{1},\ldots ,a_{d})$ is the induced
subgraph of $\DL$ on the set 
\begin{equation*}
\{ x_{1}\cdots x_{d} \in \DL :a_{j}\lle x_{j}\forall j=1,\ldots,d\} 
\end{equation*}
with \emph{height\/} $\he=\he(\Pol) $.
\end{dfn}

Note that $\Pol$ is finite; the octahedra $\Oct$ of
Definition~\ref{octahedron} are unions of $2^d$ polyhedra of the same
height. Set
\begin{equation}\label{eq:Bn}
\begin{aligned}
\hh &= \hh(\Pol) = \bigl(\hor(a_1),\ldots, \hor(a_d)\bigr) \in \Z^d\,,\\ 
B_{\hh} &= \{ \kk \in \A_{d-1} : 
k_j \ge \hor(a_j)\;\text{for}\;j =1,\ldots, d \}\,.
\end{aligned}
\end{equation}
If $\Pol$ \and $\wt \Pol$ are two polyhedra with the same height 
$\he=\he(\Pol)=\he(\wt\Pol)$, then the graphs $B_{\hh(\Pol)}$ and
$B_{\hh(\wt\Pol)}$ are isomorphic. As a representative of their 
isomorphism class, we single out the following one, by a slight abuse of
notation: 
$$
B_{\he} = B_{(0,\ldots,0,-\he)} =  \{ \kk \in \A_{d-1} : 
k_1, \ldots, k_{d-1} \ge 0\,,\; k_d \ge -\he \}\,.
$$
It contains $\binom{\he+d-1}{d-1}$ elements.
If $x=x_{1}\cdots x_{d}$ is an element of $\Pol(a_1,\dots,a_d)$ and 
$\Hor(x)=\bigl(\hor(x_1),\ldots, \hor(x_d)\bigr)$ then 
$\Hor(x) - \hh - \he\ee_d \in B_\he $. In particular, for each $j$, we have
$x_j \in T(a_j)= \{ y_j \in T_j :a_j\lle y_j \}$ and $d(x_j,a_j) \le \he\,.$ 
The \emph{boundary} $\bd \Pol$ of $\Pol$, that is, the set of all 
points in $\Pol$
having a neighbour in $\DL \setminus \Pol$, consists of all points $x \in \Pol$
for which there is at least one $j$ such that $x_j = a_j$. The \emph{interior}
of $\Pol$ is $\Pol^o = \Pol \setminus \bd \Pol$. Analogously, we set 
$\bd B_{\hh} = \{\kk \in B_{\hh} : k_j=\hor(a_j) \;\text{for some}\;j\}$ and 
$B_{\hh}^o = B_{\hh} \setminus \bd B_{\hh}$, 
the boundary and interior of the graph $B_{\hh}$, compare with
Figure~7.
$$
\beginpicture
\setcoordinatesystem units <4mm,3.5mm> 

\setplotarea x from -14 to 14, y from 0 to 14

\put {$B_{\he}$ for $d=3$} [l] at -14 6

\plot -6 0   6 0 /         
\plot -5 2   5 2 /      
\plot -4 4   4 4 /           
\plot -3 6   3 6 /      
\plot -2 8   2 8 /         
\plot -1 10  1 10 /         

\plot  -4 0  -5 2 /    
\plot  -2 0  -4 4 /     
\plot   0 0  -3 6 /      
\plot   2 0  -2 8 /     
\plot   4 0  -1 10 /
\plot   6 0   0 12 /    

\plot   4 0   5 2 /    
\plot   2 0   4 4 /     
\plot   0 0   3 6 /      
\plot  -2 0   2 8 /     
\plot  -4 0   1 10 /
\plot  -6 0   0 12 /    

\put {$(\he,0,-\he)$} [r] at -6.2 0
\put {$(0,\he,-\he)$} [l] at  6.2 0
\put {$(0,0,0)$} [c] at  0 12.7 


\endpicture
$$
\begin{center}
\centerline\emph{Figure 7}
\end{center}

\vspace{.1cm}

For $\kk \in B_{\hh}$, the $\kk$-th \emph{level} of $\Pol$ is defined by
$$
\begin{aligned}
L_{\kk} =L_{\kk}(a_{1},\ldots ,a_{d}) 
&=\{ x \in \Pol :\hor(x_{j}) = k_{j}\,,\;j=1,\ldots ,d \} \\
&=\prod_{j=1}^{d}\left( T_j(a_j) \cap 
H_{k_{j}}^{j}\right)\,
\end{aligned}
$$
where $H_{m}^{j}$ denotes the $m$-th horocycle of the tree $T_{j}$. 
As in Figure~1, we can label all edges of $T_j$ by the elements of $\Zyk_{q_j}$
such that for each vertex of $T_j$, the edges to its successors carry 
all distinct labels. Write $v_{j,l}$ for the successor of $a_{j}$ such that 
the edge from $a_{j}$ to $v_{j,l}$ has label $l\in \Zyk_{q_j}$.
We choose a function $\varphi^{j}:\Zyk_{q_j}\to\C$  with
\begin{equation}\label{eq:phijl}
\sum_{l\in \Zyk_{q_j}}\varphi^{j}(l) =0 \AND 
\sum_{l\in \Zyk_{q_j}}\left( \varphi^{j}(l) \right)^{2}=1\,.
\end{equation}
For $j \in \{1,\ldots, d\}$ and $k \in \Z$, we denote by
$f_{k}^{j}=f_{k}^{j}[a_j,\varphi^j]: T_j \to \R$ the function 
\begin{equation}\label{eq:fjk}
f_k^j(x_j) = \begin{cases} 
         \varphi^j(l) \,q^{(\hor(a_j)-k+1)/2}\,, &\text{if}\; k > \hor(a_j)
	 \;\text{and}\;v_{j,l}\lle x_j\in H_{k}^{j} \\ 
         0\,, & \text{otherwise}
              \end{cases}
\end{equation}
and we define for each $\kk = (k_1, \ldots, k_d) \in \A_{d-1}$ the
following function on $\DL\,$:  
\begin{equation}\label{eq:fkn}
f_{\kk,\Pol}(x) =f_{\kk}[\Pol,\varphi^1,\ldots ,\varphi^d](x) =
\prod_{j=1}^{d}f_{k_j}^{j}(x_{j})\,,\quad x=x_1\cdots x_d \in \DL\,.
\end{equation}
Since $f_{\kk,\Pol} \equiv 0$ on $\DL \setminus \Pol^o$, we can also consider
$f_{\kk,\Pol}$ as a function on $\Pol=\Pol(a_{1},\ldots ,a_{d})$ which is $0$
on $\bd \Pol$.

\begin{lem}\label{lem:Pfks}
The functions $f_{\kk,\Pol}$, where $\kk \in B_{\hh}^o$ and $\hh = \hh(\Pol)$, 
are  horizontal and orthonormal in $\ell^{2}(\DL)$, and
\begin{equation*}
Pf_{\kk,\Pol}=\frac{1}{D}\sum_{i=1}^{d}\sum_{j\neq i}\sqrt{q_{i}q_{j}}\,
f_{\kk +\ee _{i}-\ee _{j},\Pol}\,.
\end{equation*}
\end{lem}

\begin{proof} It is immediate by construction that the $f_{\kk,\Pol}$,
$\kk \in B_{\hh}^o$, are horizontal and orthonormal. 
Regarding the action of $P$, first note that for each $x_j \in T_j$, we have 
$$
\sum_{y_j: y_j^-=x_j} f_{k}^j(y_j)= \sqrt{q_j}\,f_{k-1}^j(x_j)\,,
$$
and when $x_j^- \ne a_j$ then also 
$f_{k}^j(x_j^-) = \sqrt{q_j}\,f_{k+1}^j(x_j)\,.$
Thus, we have for $x \in \DL \setminus \bd \Pol$
$$
\begin{aligned}
Pf_{\kk ,\Pol}(x) 
&=\frac{1}{D}\sum_{i=1}^{d}\sum_{j\neq i}\sum_{y\in N_{j,i}(x)} 
  f_{\kk,\Pol}(y) \\
&=\frac{1}{D}\sum_{i=1}^{d} f_{k_{i}}^{i}(x_i^-) \sum_{j\neq i}
  \Biggl( \sum_{\, y_j \in T_j\,:\, y_j^-=x_j} f_{k_{j}}^{j}(y_{j}) \Biggr) 
  \prod_{l\neq i,j} f_{k_{l}}^{l}(x_{l})
  \\[4pt]
&\text{(using now that $x \notin \bd\Pol$)}\\
&=\frac{1}{D} \sum_{i=1}^{d} \sqrt{q_i}\,f_{k_i+1}^{i}(x_i)
  \sum_{j\neq i}  \sqrt{q_j}\,f_{k_j-1}^j(x_j)
  \prod_{l\neq i,j} f_{k_{l}}^{l}(x_{l}) \\
&= \frac{1}{D} \sum_{i=1}^{d} \sum_{j\neq i} 
   \sqrt{q_iq_j}f_{\kk +\ee _{i}-\ee_{j},\Pol}(x)\,.
\end{aligned}
$$
When $x \in \bd \Pol$, there is $j'$ such that 
$x_{j'} = a_{j'}$. Thus, the double sum in the last line vanishes. 
On the other hand, we split the sum in the second line into three pieces. 
The first is the sum over all pairs $(i, j)$ with $i \ne j$  and both 
$i,j \ne j'$. Then $\prod_{l\neq i,j} f_{k_{l}}^{l}(x_{l}) =0$, since
this product contains the factor $f_{k_{j'}-1}^j(a_{j'}) =0$.
The second is the sum over all pairs $(i,j)$  with $i \ne j'$
and $j = j'$. Then  
$\sum_{y_j^-=x_j} f_{k_{j}}^{j}(y_{j}) = 0$ since 
$x_{j'} = a_{j'}$. The third is the sum over all pairs $(i,j)$
with $i = j'$ and $j \ne j'$. But then 
$f_{k_{i}}^{i}(x_i^-) = f_{k_{j'}}^{j'}(a_{j'}^-) =0$.
Thus, we also have $Pf_{\kk ,\Pol}(x) =0$ when $x \in \bd \Pol$.
\end{proof}

\begin{cor}\label{cor:Pfkseqq}
If $q_{1}=\ldots =q_{d}=q$ then for all $\kk \in B_{\he}^o$,
\begin{equation*}
Pf_{\kk ,\Pol}=\frac{1}{d(d-1)}\sum_{i=1}^{d}\sum_{j\neq i}
f_{\kk + \ee_{i}-\ee _{j},\Pol}\,.
\end{equation*}
\end{cor}
$\,$\\[2pt]
{\bf B. The spectra of $Q$ and $P$.}
\\[2pt]
Lemma \ref{lem:Pfks} leads us to consider the self-adjoint convolution operator
on $\ell^2(\Z^d)$ defined by
\begin{equation}\label{eq:Qdef}
Q\ff(\kk) = \frac{1}{D}\sum_{i=1}^{d}\sum_{j\neq i}\sqrt{q_{i}q_{j}}\,
\ff(\kk +\ee _{i}-\ee _{j})\,. 
\end{equation}
It leaves the subspace $\ell^2(\A_{d-1})$ invariant.

\begin{pro}
The spectrum of $Q$ is an interval $[\rho'(Q)\,,\,\rho(Q)]$ with endpoints
$\rho(Q) = \sum_{i,j: j \ne i} \sqrt{q_iq_j}/D$ and 
$\rho'(Q) \ge -1/(d-1)$.

\smallskip

{\rm (i)} If there is an $i$ such that $\sqrt{q_{i}}>\sum_{j\neq i}\sqrt{q_{j}}$
then $\rho'(Q) > -1/(d-1)$.

\smallskip
  
{\rm (ii)} If $d=2$ then $\rho'(Q)=-\rho(Q)$.

\smallskip

{\rm (iii)} If for each $i$ there is $j \ne i$ such that $q_i=q_j$ then 
$\rho'(Q) = -1/(d-1)$. 

\smallskip

{\rm (iv)} In particular, if $q_1 = \ldots = q_d$ then 
$\spec(Q)=[-\frac{1}{d-1}\,,\,1]\,$.
\end{pro}

\begin{proof}
For $\tb = (t_1, \ldots, t_d) \in [0\,,\,2\pi]^d$, the conjugate of
the Fourier transform of the operator $Q$
is the operator of multiplication with the function
$$
\wh Q(\tb) = \frac{1}{D} \sum_{i,j: j \ne i} \sqrt{q_iq_j} \,\cos(t_i-t_j)\,.
$$
It is well known that $\spec(Q)$ coincides with the set of values of
$\wh Q$. As a continuous image of the connected set $[0\,,\,2\pi]^d$,
it must be an interval.
The upper bound $\rho(Q)$ is attained for $\tb=\zero$. 
For the lower bound, rewrite
$$
\wh Q(\tb) = \frac{1}{D}
\left( \left| 
 \sum_{j=1}^{d}\sqrt{q_{j}}\,e^{\im t_{j}}\right|^{_{\displaystyle 2}}
- \sum_{j=1}^{d}q_{j}\right) 
\ge -\frac{1}{D}\sum_{j=1}^{d}q_{j}=-\frac{1}{d-1}\,.
$$
(Here, $\im$ is the complex unit.)

In case (i), we have
$$
\left| \sum_{j=1}^{d}\sqrt{q_{j}}\,e^{\im t_{j}}\right|
\geq \sqrt{q_{i}}-\sum_{j\neq i}\sqrt{q_{j}}>0\,,
$$
whence $\inf_{\tb} \wh Q(\tb) > -1/(d-1)$.  

In case (ii), we find $\wh Q(\pi,0) = -\rho(Q)$.

In case (iii), suppose without loss of generality that there are
$2 \le r(1) \le \ldots \le r(s) = d$ such that $r(l)-r(l-1) \ge 2$
and $q_j = q_{r(l)}$ for $r(l-1) < j \le r(l)$. For those $j$, we
set $t_j = 2\pi j/(r(l)-r(l-1))$, so that $\sum_{j=r(l-1)+1}^{r(l)}e^{\im t_{j}} = 0$,
whence $\wh Q(t_1,\ldots,t_d) = -1/(d-1)$.
\end{proof}

Now consider the vector space $\Vf_0(B_{\hh})$ of all functions 
$\ff: B_{\hh} \to \C$
with $\ff \equiv 0$ on $\bd B_{\hh}$, and the symmetric operator (matrix)
$Q_{\hh}: \Vf_0(B_{\hh}) \to \Vf_0(B_{\hh})$ defined by 
$$
Q_{\hh}\ff(\kk) = \begin{cases} Q\ff(\kk)\,,&\text{if}\;\kk \in B_{\hh}^o\,,\\
                            0\,,&\text{if}\;\kk \in \bd B_{\hh}\,.
 	      \end{cases}
$$	      	
Again by slight abuse of notation, we write $Q_{\he}$ for $Q_{(0,\dots,0,-\he)}$,
acting on $\Vf_0(B_{\he})$. Since $B_{\hh}$ and $B_{\he}$ are isomorphic when
$\langle \hh, \uno \rangle=-\he$, we have $\spec(Q_{\hh}) = \spec(Q_{\he})$ in
this case.
 
\begin{lem}\label{lem:specunion} \hspace{1.3cm}
$
\displaystyle 
\spec(Q) = \overline{\bigcup_{\he  \ge 2} \spec(Q_{\he})}.
$  
\end{lem}

\begin{proof} The operator $Q$ is the simple random walk operator on $\A_{d-1}$.  In particular, 
$\spec(Q)$ is the spectrum of $Q$ acting on $\ell^2(\A_{d-1})$.

Consider functions in  $\Vf_0(B_{\he})$ as functions on $\A_{d-1}$ with value
$0$ outside of $B_{\he}^o$. For each $\ff \in  \Vf_0(B_{\he})$, we have 
$Q_{\he}\ff(\kk) = Q\ff(\kk)$ for each 
$\kk \in \A_{d-1} \setminus \bd B_{\he}$, while $\ff \equiv 0$ on $\bd B_{\he}$.
Therefore 
$$
\langle Q_{\he} \ff, \ff \rangle = \langle Q \ff, \ff \rangle 
\begin{cases} \le \rho(Q) \langle \ff, \ff \rangle & \\
               \ge \rho'(Q) \langle \ff, \ff \rangle\,. &  
\end{cases}	       
$$ 
Consequently $\rho'(Q) \le \rho'(Q_{\he}) \le \rho(Q_{\he}) \le \rho(Q)$,
and 
$$
\spec(Q_{\he}) \subset [\rho'(Q)\,,\,\rho(Q)] = \spec(Q)\,.
$$
On the other hand, choose $\lambda \in \spec{Q}$ and $\varepsilon > 0$. 
As $Q$ is self-adjoint there exists a \emph{finitely supported} 
$\ff \in \ell ^{2}\left( \A _{d-1}\right) $ of norm~$1$ such that
$\left\| \lambda \ff -Q\ff \right\| <\varepsilon$. There is
an $\hh $ such that $\ff\in \Vf_{0}(B_{\hh})$ and $Q\ff = Q_{\hh}\ff$. Thus,
expanding $\ff$ with respect to an orthonormal basis of $\Vf_{0}(B_{\hh})$
consisting of $Q_{\hh}$-eigenfunctions, we see that there must be
an eigenvalue $\la'$ of $Q_{\hh}$ such that $|\la-\la'| < \varepsilon$. 
Since $\spec(Q_{\hh}) = \spec(Q_{\he})$, where 
${\he}=-\left\langle \hh, \uno \right\rangle$, the proof is complete.
\end{proof}

\begin{lem}\label{lem:perp}
Let $\Pol=\Pol(a_{1},\ldots ,a_{d})$ and 
$\wt \Pol = \Pol(\tilde{a}_{1},\ldots ,\tilde{a}_{d})$ be two polyhedra, both 
of height $\geq 2$, and set $\hh = \hh(\Pol)$, $\wt\hh=\hh(\wt\Pol)$. 
Furthermore let $\varphi^{1},\ldots ,\varphi^{d}$
and $\wt\varphi^{1},\ldots ,\wt\varphi^{d}$ be functions 
on $\Zyk_{q_j}$ satisfying \eqref{eq:phijl}. 
For $\kk \in B_{\hh}^o$ and $\lk \in B_{\wt \hh}^o$, write 
$f_{\kk,\Pol}=f_{\kk}[\Pol, \varphi^1,\ldots,\varphi^d]$ and 
$\wt{f}_{\lk,\wt\Pol}=f_{\lk}[\wt\Pol,\wt\varphi^{1},\ldots,\wt\varphi^{d}]\,$. 

Suppose that one of the following conditions is satisfied.

\begin{enumerate}
\item  $(a_{1},\ldots ,a_{d}) \neq (\tilde{a}_{1},\ldots ,\tilde{a}_{d})\,,$
or
\item  $(a_{1},\ldots ,a_{d}) = (\tilde{a}_{1},\ldots ,\tilde{a}_{d})$ and 
$\varphi^{j}\perp \wt\varphi^{j}$ in $\ell^2(\Zyk_{q_j})$ for some
$j \in \{1,\ldots ,d\}$.
\end{enumerate}

Then 
\begin{equation*}
\{ f_{\kk,\Pol}: \kk \in B_{\hh}^o\} \perp 
\{ \wt{f}_{\lk ,\wt \Pol} : \lk \in B_{\wt\hh}^o\}\,. 
\end{equation*}
\end{lem}

\begin{proof}
(1)  If $\Pol \cap \wt \Pol =\emptyset$ then the above two sets are
obviously orthogonal.

Suppose now $\Pol \cap\wt \Pol \neq \emptyset$. There is $j$ such that $a_j
\ne \tilde a_j\,$, and $a_j,\tilde a_j$ are comparable in the ancestor
relation $\lle$. Without loss of generality, suppose $a_1\neq\tilde
a_1$ and $a_{1}\lle\tilde{a}_{1}$ and $a_{1}\neq \tilde{a}_{1}$. Then
$\tilde{a}_{1} \succcurlyeq v_{1,l}$, where the latter is one of the
successors of $a_1$ in $T_1\,$.

Let $\kk\in B_{\he}^o$ and $\lk \in B_{\tilde {\he}}^o$ be given. If
$\kk\neq \lk$ then certainly $f_{\kk ,\Pol}\perp \wt f_{\lk,\wt\Pol}$,
as these functions are supported on disjoint sets.

Assume therefore $\kk=\lk$.
By construction \eqref{eq:fjk}, the function $f_{k_{1}}^{1}$
on the subtree of $T_1$ is constant on the support of 
$\wt f_{k_{1}}^{1}$. Also by construction, the sum of 
$\wt f_{k_{1}}^{1}$ over its support is $=0$. 
Thus, 
$$
\sum_{x_1 \in T_1: \hor(x_1) = k_1} 
f_{k_{1}}^{1}(x_1)\,\wt f_{k_{1}}^{1}(x_1) = 0\,,
$$
that is, $f_{k_{1}}^{1}\perp \wt f_{k_{1}}^{1}$. 
Hence, again by construction \eqref{eq:fjk},
\begin{equation*}
\sum_{x}f_{\kk ,\Pol}(x) \wt f_{\lk,\wt\Pol}(x) =
\prod_{j}\sum_{x_{j}}f_{k_{j}}^{j}(x_{j}) \wt f_{k_{j}}^{j}(x_{j}) =0\,.
\end{equation*}

(2) If $(a_{1},\ldots ,a_{d}) = (\tilde{a}_{1},\ldots ,\tilde{a}_{d})$ 
then $\wt \Pol = \Pol$ and $\wt\hh=\hh$. If $\kk, \lk \in B_{\he}^o$
are distinct then $f_{\kk,\Pol}$ and $\wt f_{\lk,\Pol}$ have disjoint 
support and are perpendicular. In the last remaining case, we compute 
\begin{equation*}
\bigl\langle f_{\kk ,\Pol},\wt f_{\kk ,\Pol}\bigr\rangle
=\prod_{j=1}^{d}\bigl\langle \varphi ^{j},\wt\varphi^{j}\bigr\rangle =0\,.
\qedhere
\end{equation*}
\end{proof}

Now let 
$\{ \psi_{\mm,{\he}} : \mm \in B_{\he}^o \}$
be an orthonormal basis of $\Vf_0(B_n)$ consisting of
eigenfunctions (-vectors) of $Q_{\he}$ with associated eigenvalues
$\la_{\mm,{\he}}\,$, parametrized by all $\mm \in B_{\he}^o$. The following is
immediate from Lemma \ref{lem:Pfks} and the fact that for 
$\langle \hh,\uno \rangle = -{\he}$, the natural
isomorphism $B_{\he} \to B_{\hh}$  is given by 
$\kk \mapsto \kk + \hh + {\he}\ee_d\,$.

\begin{cor}\label{cor:gmn} Let $\Pol=\Pol(a_1,\ldots,a_d)$, $\hh=\hh(\Pol)$ 
and ${\he}={\he}(\Pol)$.
The functions $g_{\mm,\Pol}:\DL \to \R$, defined by
\begin{equation*}
g_{\mm,\Pol}=\sum_{\kk\in B_{\he}^o}
\psi _{\mm,{\he}}(\kk)\, f_{\kk+\hh+{\he}\ee_d,\Pol}\,,\quad
\mm \in B_{\he}^o\,,
\end{equation*}
are orthonormal and horizontal vectors in $\ell^{2}(\Pol)$ as well as
$\ell ^{2}(\DL)\,$. They satisfy 
\begin{eqnarray*}
\spn\left\{ g_{\mm,\Pol}:\mm\in B_{\he}^o\right\}  &=&
\spn\left\{ f_{\kk ,\Pol}:\kk\in B_{\hh}^o\right\} \AND\\
Pg_{\mm,\Pol} &=&\la_{\mm,{\he}}\cdot g_{\mm,\Pol}.
\end{eqnarray*}
\end{cor}

Remember that each $g_{\mm,\Pol}$, $\mm \in B_{\he}^o\,$, depends 
not only on the polyhedron $\Pol=\Pol(a_{1},\ldots ,a_{d})$, but also
on the functions $\varphi_{1},\ldots ,\varphi_{d}$ of \eqref{eq:phijl}, 
that is, 
\begin{equation}\label{eq:gdepends}
g_{\mm,\Pol}=g_{\mm}[\Pol,\varphi_{1},\ldots ,\varphi_{d}]\,.
\end{equation}

For each $q_j$, we now select functions $\varphi^j_k$, 
$k \in \Zyk_{q_j} \setminus \{0\}$ that satisfy \eqref{eq:phijl} and
are mutually orthogonal:
\begin{equation}\label{eq:phichoice}
\varphi^j_l(s) = \begin{cases} 0\,,&s = 0, \dots, l-2\,,\\
       (q_j-l)\big/\sqrt{(q_j-l)(q_j+1-l)}\,,&s=l-1\,,\\[3pt]
       -1\big/\sqrt{(q_j-l)(q_j+1-l)}\,,&s=l, \dots, q_j-1\,.
        \end{cases}
\end{equation}

\begin{pro}\label{pro:polbasis} The set
\begin{equation*}
\mathfrak{B}_{\Pol} =
\bigl\{ g_{\mm}[\wt\Pol,\varphi^1_{l_1},\ldots,\varphi^d_{l_d}] :
\wt\Pol \subset \Pol\,,\; \mm\in B_{{\he}(\wt\Pol)}^o\,,\;
l_{j}\in \Zyk_{q_j} \setminus \{0\}\bigr\}\,,
\end{equation*}
where $\wt\Pol$ runs through all polyhedra with height ${\he}(\wt\Pol)\ge 2$ 
contained in $\Pol$, is an orthonormal basis of the vector space 
$\mathfrak{V}_0{\Pol}$ of all 
functions which are horizontal and supported by $\Pol$. 
\end{pro}

\begin{proof} Lemma \ref{lem:perp} and Corollary \ref{cor:gmn}
imply that the elements of $\mathfrak{B}_{S}$ are orthonormal, whence
linearly independent. What remains is to show that $\mathfrak{B}_{S}$ 
spans the whole of $\mathfrak{V}_0(\Pol)$.
We can replace each
$g_{\mm}[\wt\Pol,\varphi^1_{l_1},\ldots,\varphi^d_{l_d}] \in \mathfrak{B}_{S}$ 
with 
$f_{\mm+\hh+{\he}\ee_d}[\wt\Pol,\varphi^1_{l_1},\ldots,\varphi^d_{l_d}]\,$, 
since the latter functions are also linearly independent and span
the same vector space by Corollary~\ref{cor:gmn}.

Recall that $L_{\kk}$ denotes the $\kk$-th level of the polyhedron $\Pol$
with height ${\he}={\he}(\Pol)$. A function $f$ in the space of all 
complex functions supported by $L_{\kk}$ is horizontal if and only if it
satisfies each of the equations
$$
\sum_{y_j \in \T_j\,: \,\hor(y_j)=\hor(x_j)}  
f(x_1\cdots x_{j-1}y_jx_{j+1}\cdots x_d) = 0\,, 
$$
for all $j \in \{1, \ldots, d\}$ and $x_i \in T_i(a_i) \cap H^j_{k_i}$
for each $i \ne j$. The number of these equations is 
$\sum_{j=1}^{d} \prod_{i\ne j} q_{i}^{k_{i}-\hor(a_i)}$, but they are
not independent. By inclusion-exclusion, we find
$$
\dim \mathfrak{V}_0(L_{\kk})
= \prod_{i=1}^{d}\bigl(q_i^{k_i-\hor(a_i)}-1\bigr)\,.
$$
The proof will be complete when we can show that 
$\mathfrak{V}_0(L_{\kk})$ is spanned by all the functions
$f_{\kk}[\wt\Pol,\varphi^1_{l_1},\ldots,\varphi^d_{l_d}]\,$, where 
$\wt\Pol \subset \Pol$, which are supported in $L_{\kk}$. To this
end, we count all functions of the latter type (which are linearly
independent). 

First of all, the interior of $\wt\Pol = \Pol(\tilde a_1,\ldots,\tilde a_d)$ 
has to intersect $L_{\kk}$. This means that for each $j$, we need to have 
$\hor(a_j) \le \hor(\tilde a_j) < k_j$. There are
$\sum_{r=0}^{k_j-\hor(a_j)-1} q_j^r$ points $\tilde a_j$ of this type. 
Thus, the number of feasible polyhedra $\wt\Pol$ is 
$
\prod_{j=1}^{d}\bigl(q_j^{k_j-\hor(a_j)}-1\bigr)/(q_j-1)\,.
$

The remaining choices that we have are those of the functions
$\varphi^j_l$, where $l \in \Zyk_{q_j} \setminus \{0\}$ and
$j\in\{1,\ldots,d\}$.  There are $\prod_{j=1}^{d}(q_j-1)$ such
choices, so that the total number of functions
$f_{\mm}[\wt\Pol,\varphi^1_{l_1},\ldots,\varphi^d_{l_d}]$ supported in
$L_{\kk}$ is $\prod_{j=1}^{d}\bigl(q_j^{k_j-\hor(a_j)}-1\bigr)$.  This
number coincides with the dimension of $\mathfrak{V}_0(L_{\kk})$, as
claimed.
\end{proof}

We now can show that as sets, the spectra of $P$ on $\ell^2(\DL)$ and
$Q$ on $\A_{d-1}$ coincide, although their ``inner structure'' is
completely different. Indeed, $\spec(Q)$ is continuous, i.e.,
the associated spectral (Plancherel) measure is absolutely continuous with
respect to Lebesgue measure, while for $\spec(P)$, we have the following.

\begin{thm}\label{thm:spectrum} 
The spectrum of the operator $P$ is pure point. It is the closure of
the set of eigenvalues
$$
\spec_p(P)= \Bigl\{ \la_{\mm,{\he}} : {\he} \ge 2\,, \mm \in B_{\he}^o \Bigr\} 
= \bigcup_{{\he} \ge 2} \spec(Q_{\he})\,.
$$
Each eigenvalue has infinite multiplicity.
An associated orthonormal basis consisting of  
finitely supported eigenfunctions of 
$P$ is given by
\begin{equation*}
\mathfrak{B} =\Bigl\{ g_{\mm}[\Pol,\varphi^1_{l_1},\ldots,\varphi^d_{l_d}] :
\Pol \;\text{polyhedron in}\; \DL\;\text{with}\;n(\Pol) \ge 2\,,\; 
\mm\in B_{{\he}(\Pol )}^o\,,\; l_{j}\in \Zyk_{q_j} \setminus \{0\} \Bigr\}\,,
\end{equation*}
where $g_{\mm}[\Pol,\varphi^1_{l_1},\ldots,\varphi^d_{l_d}]$ is 
defined in Corollary~\ref{cor:gmn} and (\ref{eq:gdepends}). 
\end{thm}

\begin{proof}
Pick any function $f$ on $\DL$ which is horizontal. Its support is then
contained in an appropriate polyhedron $\Pol$. Proposition~\ref{pro:polbasis}
implies that $f$ is in the span of $\mathfrak{B}_{\Pol}$. Hence every 
horizontal function is a linear combination of functions in $\mathfrak{B}$.
Therefore $\mathfrak{B}$ is an orthonormal basis by Lemma~\ref{lem:dense}, 
and the whole spectrum must consist of the closure of the set of associated 
eigenvalues.
\end{proof}

In particular, we get $\rho(P)=\rho(Q)$, that is, the spectral radius of $P$
is
\begin{equation}\label{eq:specrad}
\rho(P) = \frac{2}{D} \sum_{i,j: i < j} \sqrt{q_iq_j}\,,
\quad\text{where}\quad D=(d-1)\sum_i q_i\,.
\end{equation}
This formula might also have been derived in a completely different
way, by use of Theorem~1(b) of~\cite{SoWo}. 

For $\DL(q_1,q_2)$, the eigenvalues of $P$ can be computed
explicitly~\cite{BaWo}. In general, this amounts to the explicit computation
of the eigenvalues of $Q_{\he}$ on $B_{\he}$. For $d \ge 3$, this is the precise
discrete analogue to computing the eigenvalues of the Laplacian on the 
$(d-1)$-dimensional Euclidean simplex whose side lengths are proportional
to $\sqrt{q_iq_j}$ (respectively). An explicit solution to the latter problem 
is known only for equilateral triangles in dimension $d-1=2$, but in no other
case: this solution goes back to the work of Lam\'e in the 19th century,
see e.g. {\sc Lam\'e~\cite{La},} and has reappeared in the literature
several times (with or without knowledge of Lam\'e's work), see  
{\sc Pinsky~\cite{Pi}, Pr\'ager~\cite{Pr1,Pr2}, McCartin~\cite{Mc}.}  
\\[6pt]
{\bf C. The spectrum of $P$ on $\DL(q,q,q)$.}
\\[2pt]
We now consider briefly  the case of $\DL(q,q,q)$, where
we have from Corollary~\ref{cor:Pfkseqq} that $Q$ is the transition operator 
of SRW on $\A_2$. 

Let $\SF_3$ be the permutation group of $\{1,2,3\}$. For $\si \in \SF_3$ and
$\kk=(k_1,k_2,k_3) \in \Z^3$, we write 
$\si\kk = (k_{\si^{-1}(1)},k_{\si^{-1}(2)}, k_{\si^{-1}(3)})$. 
Also, we write $\Alt_3$ for the subgroup of even permutations in $\SF_3$.
Now we define the following functions for ${\he} \ge 2$.
\begin{equation}\label{eq:psimndef}
\psi_{\mm,{\he}}(\kk) = \frac{1}{\sqrt{6}\,{\he}}\sum_{\si \in \SF_3}
\func{sign}(\si) \exp \left( \func{sign}(\si) \frac{2\pi \im}{3{\he}}
\langle \mm,\si \kk  \rangle \right)\,,\quad \mm\,,\;\kk \in \A_2\,.
\end{equation}
\begin{pro}\label{pro:psimnprop} We have the following.
\begin{align}
\psi_{\mm,{\he}} &= 0 
\quad \text{on $\;\bd B_{\he}\;$ for each $\;\mm \in B_{\he}^o\,$} \tag{a}\\
\langle \psi_{\mm,{\he}}\,,\,\psi_{\mm^{\prime},{\he}}\rangle &= 
\delta_{\mm}(\mm^{\prime}) 
\quad \text{on $\;B_{\he}\;$ for all 
$\;\mm\,,\;\mm^{\prime}\in B_{\he}^o\,$,} 
\tag{b} \\
Q_{\he}\,\psi_{\mm,{\he}} &= \la_{\mm,{\he}} \cdot \psi_{\mm,{\he}} \quad 
\text{on $\;B_{\he}\,$, where} \tag{c}
\end{align}
$$
\la_{\mm,{\he}}= \frac{1}{3}\left( 
\cos \biggl( \frac{2\pi (m_{1}-m_{2})}{3{\he}}\biggr) + 
\cos \biggl( \frac{2\pi (m_{2}-m_{3})}{3{\he}}\biggr) + 
\cos \biggl( \frac{2\pi (m_{3}-m_{1})}{3{\he}}\biggr) 
                         \right).
$$
\end{pro}

\begin{proof} (a) The boundary of $B_{\he}$ consists of the points
$k(\ee_1-\ee_3)$ (first part), $k(\ee_2-\ee_3)$ (second part) 
and ${\he}(\ee_1-\ee_3)+k(\ee_2-\ee_1)=({\he}-k,k,-{\he})$ (third part),
where $k=0,\ldots, \he $. 

Let $\tau_{i,j} \in \SF_3$ be the transposition of $i$ and $j$ ($i\ne j$).
Then 
$$
\begin{aligned} &\psi_{\mm,{\he}}(k\ee_i-k\ee_j)\\
&= \frac{1}{\sqrt{6}\,{\he}}\sum_{\si \in \Alt _3} \left(
\exp\left(\frac{2\pi \im}{3{\he}} \langle \mm,\si(k\ee_i-k\ee_j) \rangle\right) - 
\exp\left(\frac{2\pi \im}{3{\he}} 
                 \langle \mm,-\si \tau_{i,j}(k\ee_i-k\ee_j)\rangle \right)
                                         \right) = 0		 
\end{aligned}
$$
Thus, $\psi_{\mm,{\he}} = 0$ on the first and second parts of $\bd B_{\he}\,$. 
Regarding the third part, note that for $\mm \in \A_2\,$,
$$
\bigl\langle \mm, ({\he}-k,k,-{\he}) \bigr\rangle +
\bigl\langle \mm, \tau_{1,2}({\he}-k,k,-{\he}) \bigr\rangle =
(m_1+m_2-2m_3){\he} \equiv 0 \,\func{mod}\, 3{\he}\,.
$$
Therefore, pairing each element $\si \in \Alt_3$ with $\si \tau_{1,2}$ as in
the above sum
yields that $\psi_{\mm,{\he}} = 0$ on the third part of $\bd B_{\he}\,$.
 
\smallskip

(b) We have by construction 
$$
\psi_{\mm,{\he}}(\kk) 
= \func{sign}(\si) \psi_{\mm,{\he}}\bigl(\func{sign}(\si) \si\kk\bigr)\,,
$$
and since for $\mm \in \A_2$, each of $\langle \mm, ({\he},{\he},-2{\he}) \rangle\,$,
$\langle \mm, ({\he},-2{\he},{\he}) \rangle$ and 
$\langle \mm, (-2{\he},{\he},{\he}) \rangle$
is a multiple of $3{\he}$, we get 
$$
\psi_{\mm,{\he}}(\kk) =\psi_{\mm,{\he}}\bigl(\kk+({\he},{\he},-2{\he})\bigr) 
=\psi_{\mm,{\he}}\bigl(\kk+({\he},-2{\he},{\he})\bigr) 
=\psi_{\mm,{\he}}\bigl(\kk+(-2{\he},{\he},{\he})\bigr).
$$
That is, up to a change of the sign, $\psi_{\mm,{\he}}$ is invariant under
the group of motions of $\A_2$ generated by the translations 
$\kk \mapsto \kk+{\he}(\ee_i+\ee_j-2\ee_l)\,$, where $\{i,j,l\} = \{1,2,3\}$,
and the reflections 
$\kk \mapsto - \tau_{i,j}\kk$, where $1 \le i < j \le 3$.
Now $B_{3{\he}}$ is the union of $9$ distinct images of $B_{\he}$ under elements
of that group, and those images meet only at their respective boundaries,
where $\psi_{\mm,{\he}}=0$. Therefore 
$$
\langle \psi_{\mm,{\he}}, \psi_{\mm^{\prime},{\he}} \rangle_{B_{3{\he}}} = 
9\,\langle \psi_{\mm,{\he}}, \psi_{\mm^{\prime},{\he}} \rangle_{B_{\he}}\,,
$$
where the inner products are taken over $B_{3{\he}}$ and $B_{\he}$, respectively.
Also, $\psi_{\mm,{\he}}=0$ on $\bd B_{3{\he}}$. 
Recall that $\sum_{k=0}^{3{\he}-1} \exp \bigl(\frac{2\pi \im}{3{\he}}rk\bigr) 
= 0$ when $r \in \Z$ is not a multiple of $3{\he}$. 
We now obtain, setting 
$\kk = (k_1,k_2,-k_1-k_2)$ if $k_1$ and $k_2$ are given,
$$
\begin{aligned}
\langle \psi_{\mm,{\he}}, &\psi_{\mm^{\prime},{\he}} \rangle_{B_{\he}} \\
&= \frac{1}{54{\he}^2} 
\sum_{\si ,\tau \in \SF }\func{sign}(\si\tau) \!\!\!\sum_{k_{1},k_{2}=0}^{3{\he}-1}
\exp \!\left(\frac{2\pi \im}{3{\he}}
\Bigl\langle \func{sign}(\si) \sigma ^{-1}\mm 
            -\func{sign}(\tau) \tau ^{-1}\mm^{\prime}\,,\,\kk \Bigr\rangle
     \!\right)  \\
&=\frac{1}{6}\sum_{\si ,\tau \in \SF}
\func{sign}(\si^{-1}\tau)
\de_{\mm}\bigl(\func{sign}(\sigma^{-1} \tau) \si\tau^{-1}\mm^{\prime}\bigr)\\
&=\sum_{\si \in \SF }\func{sign}(\si)
\de_{\mm}\bigl((\func{sign}(\si) \si \mm^{\prime}\bigr)
=\de_{\mm}(\mm^{\prime})\,,
\end{aligned}
$$
as $B_{\he}^o \cap \bigl(\func{sign}(\si) \si\bigr) B_{\he}^o = \emptyset$ 
for $\si \neq \id$.

\smallskip

(c) is a straightforward computation.
\end{proof}

\begin{cor}\label{cor:DL3spec} The spectrum of the SRW operator $P$ on
$\DL(q,q,q)$ is the interval $\bigl[ -\frac12\,,\,1\bigr]$. It is the closure
of the set of eigenvalues 
$$
\{ \la_{\mm,{\he}} : {\he} \ge 2\,,\; \mm \in \A_2\;\text{with}\;\, 
m_1, m_2 > 0\,,\; m_3 > -{\he} \}
$$
given by Proposition~\ref{pro:psimnprop}(c).
\end{cor}

The associated orthonormal basis of finitely supported eigenfunctions can
be computed explicitly by use of the functions $\psi_{\mm,{\he}}$ of
\eqref{eq:psimndef}, according to Theorem~\ref{thm:spectrum}.

We remark here that the groups $\Ga_3(\Zq)$ and 
$\Ga_d(\F_q)$ ($d \ge 3$) constructed in 
Corollaries \ref{cor:3trees} and \ref{cor:dtreesCayley} are the first examples 
of finitely \emph{presented} groups where SRW has a pure point spectrum. 
\\[6pt]
{\bf D. Spectral measure and return probabilities on $\DL(q,q,q)$.}
\\[2pt]
If $\mathfrak{B}$ is the orthonormal basis of Theorem~\ref{thm:spectrum}
and we write $\la(g)$ for the eigenvalue associated
with $g \in \mathfrak{B}$, then we see that the $n$-step transition
probabilities $p^{(n)}(x,y) = p^{(n)}(y,x) =\langle P^n\de_x,\de_y\rangle$
satisfy
\begin{equation}\label{eq:specsum}
\begin{aligned}
p^{(n)}(x,y) &= \sum_{g \in \mathfrak{B}} \la(g)^n \, g(x)g(y)
  = \sum_{\la \in \spec_p(P)} \la^n \,\mu_{x,y}(\la)\,,\quad\text{where}\\            
\mu_{x,y}(\la) &= \sum_{g \in \mathfrak{B}: \la(g)=\la} g(x)g(y)
\end{aligned}
\end{equation}
Thus, we obtain the $(x,y)$-element $\mu_{x,y}$  of the \emph{spectral measure}
of the self-adjoint operator $P$, see e.g. {\sc Grigorchuk and
\.Zuk~\cite{GrZu2}} for spectra of Markov operators on graphs, and
\cite{BaWo} in the context of $\DL$-graphs. Since $\spec(P)$ is pure point,
$\mu_{x,y}$ is an infinite sum of weighted point masses 
$\mu_{x,y}(\la)\,\de_{\la}$. We want to have a closer look at the measure
$\mu_{x,x} = \mu_{o,o}$, which plays the role of the 
\emph{Plancherel measure}.  

The eigenvalues found in Proposition \ref{pro:psimnprop} satisfy
the obvious relation $\la_{\mm,{\he}} = \la_{\ell\mm,\ell {\he}}$ for each
$\ell \in \N$. Besides, we also have $\la_{\mm,{\he}} = \la_{\tau_{1,2}\mm,{\he}}$
(exchanging of $m_1$ and $m_2$), and when ${\he}$ and $\mm \in B_{\he}^o$ 
are such that $2m_1+4m_2 = 3{\he}$ then $\la_{\mm,{\he}} = -1/3$. In general, 
it is a very delicate task to determine all $(\mm',\he')$ for which 
$\la_{\mm',\he'} = \la_{\mm,{\he}}$,
where $(\mm,{\he})$ is given. Thus, instead of a complete description, we
pull back $\mu_{o,o}$ to a measure $\wt \mu$ on the parameter space
\begin{equation}\label{eq:MC}
\MC = \{(\mm,{\he}) : {\he} \ge 2\,,\;\mm \in B_{\he}^o\,,\; 
\gcd({\he},m_1,m_2) = 1 \},
\end{equation}
by considering only the first, obvious relation above.
(Since $m_3 = -m_2-m_1$, the parametrization is in reality over all
$\he  \ge 2$ and $m_1,m_2 \ge 1$ with $m_1+m_2 \le \he -1$.) 
We note that for any polyhedron $\Pol$, we have 
$g_{\mm}[\wt\Pol,\varphi^1_{l_1},\varphi^2_{l_2},\varphi^3_{l_3}](o) \ne 0$ if 
and only if $o \in S^o$ and (since the edge from each $o_j$ to its predecessor
in $T_j$ has label $0$) if $l_1= l_2 = l_3 = 1$. We write 
$g_{\mm}[\wt\Pol] = g_{\mm}[\wt\Pol,\varphi^1_{1},\varphi^2_{1},\varphi^3_{1}]$.
Thus, for $(\mm,{\he}) \in \MC$, we compute
$$
\wt \mu(\mm,{\he}) 
= \sum_{\ell = 1}^{\infty} \SUM(\ell \mm,\ell {\he})\,,\quad \text{where}
\quad \SUM(\mm,{\he}) = \sum_{\Pol: o \in \Pol^o, {\he}(\Pol)={\he}}
\bigl(g_{\mm}[\wt\Pol](o)\bigr)^2.
$$
Now let $\Pol = \Pol(a_1,a_2,a_3)$, $o \in \Pol^o$, ${\he}(\Pol)={\he}$ and 
$\hh = \hh(\Pol)$. Then $a_j \lle o_j$, $a_j \ne o_j$ ($j=1,2,3$), and
$\kk = -\hh-{\he}\ee_3 \in B_{\he}^o$. We get from Corollary \ref{cor:gmn} that
$$
g_{\mm}[\wt\Pol](o) = \psi_{\mm,{\he}}(\kk) \, f_{\zero}(o) 
= \psi_{\mm,{\he}}(\kk)\, (q-1)^{3/2}q^{-{\he}/2}\,.
$$
If we vary $\Pol$ such that ${\he}(\Pol) = {\he}$ and $o \in \Pol^o$, then this 
amounts to varying $\kk = -\hh(\Pol) - {\he}\ee_3$ in $B_{\he}^o$. Thus
$$
\SUM(\mm,{\he}) = \sum_{\kk \in B_{\he}^o} 
\psi_{\mm,{\he}}(\kk)^2\,(q-1)^3q^{-{\he}} 
= (q-1)^3q^{-{\he}} \,,
$$
and therefore
\begin{equation}\label{eq:planch}
\wt \mu(\mm,{\he}) = (q-1)^3/(q^{\he}-1)\,.
\end{equation}
The Plancherel measure $\mu_{o,o}$ is the image of $\wt\mu$ under the
mapping $\MC \to \spec(P)$, $(\mm,{\he}) \mapsto \la_{\mm,{\he}}$. 
(If $m_1 \ne m_2$ then this mapping is at least two-to-one.) We get in 
particular the 
following nice expression for the $n$-step return probabilities.
\begin{cor}\label{cor:returnprob} With $\MC$ and $\la_{\mm,{\he}}$
defined as in \eqref{eq:MC} and Proposition \ref{pro:psimnprop}, 
respectively, the $n$-step return probabilities for SRW on $\DL(q,q,q)$
are given by
$$
\begin{aligned}
p^{(n)}(x,x) &= \sum_{(\mm,{\he}) \in \MC} 
\la_{\mm,{\he}}^n \,(q-1)^3/(q^{\he}-1)\\
&\sim A \,n^{1/6}\, \exp(-B\,n^{1/3})\,,\quad\text{as}\; n\to\infty\,,
\quad\text{for some}\ A, B > 0\,.
\end{aligned}
$$
Here $a_n\sim b_n$ means that $a_n/b_n\to1$ as $n\to\infty$.
\end{cor} 
 
We omit the computational details of the asymptotic formula,
which follow the method of \cite[\S 5]{BaWo}. Since for $d \ge 4$,
the precise information about the spectral measure is not available,
we cannot extend the above corollary to all $\DL$-graphs. However,
a rougher estimate is available.  
If $(a_n)$ and $(b_n)$ are two sequences of positive numbers, then
we say that $a_n \lle b_n$ if there are constants $C, D > 0$ such that
$a_n \le C b_{Dn}$. If $a_n \lle b_n$ and $b_n \lle a_n$, then we say
that $(a_n)$ and $(b_n)$ have the same \emph{asymptotic type,} and write
$a_n \approx b_n$. Thus, SRW on $\DL(q,q,q)$ satisfies
$p^{(n)}(x,x) \approx \exp(-n^{1/3})$ by the above, and the same is
true for $\DL(q,q)$ (when $n$ is even, since $\DL(q,q)$ is bipartite),
see \cite{Ber,BaWo}. Also, for arbitrary $d$, when
$q_1,\dots,q_d$ do not all coincide then $\rho(P) < 1$, that is,
$$
p^{(n)}(x,x) \approx \exp(-n) \quad\text{on}\quad \DL(q_1,\dots,q_d)\,,
$$
again considering only even $n$ when $d=2$.

\begin{pro}\label{pro:type} For each $d \ge 3$, SRW on $\DL_d(q)$ satisfies
$p^{(n)}(x,x) \approx \exp(-n^{1/3})$. 
\end{pro}

\begin{proof}We follow a well-known method, see e.g. \cite[\S 14--15]{Wbook}.
It is a general fact that on a vertex-transitive graph with exponential
growth, $p^{(n)}(x,x) \lle \exp(-n^{1/3})$, see 
\cite[Corollary 4.15(b)]{Wbook}. Thus, we only need to prove the lower bound.
This is based on the following estimate, where we take $x=o$, and
$A_m$ is a suitably chosen finite set of vertices of our graph, see
e.g. \cite[\S 15]{Wbook}.
$$
p^{(2n)}(o,o) \ge \Prob_o[Z_n \in A_m]^2 \big/ |A_m|\,,
$$
where $Z_n$ is the random position of the random walk starting at $o$.
(For details of the probabilistic notions, see \S \ref{sect:Poisson} below.)

Now let $a_j = o_j^{-m}$ be the $m$-th predecessor of $o_j$, and set
$A_m = \Pol(a_1, \dots,a_m)$. The cardinality of this set is of the order 
of $C\,m^{d-1}q^{dm}$. Let 
$Z_{j,n} \in T_j$ be the $j$-th coordinate of $Z_n$. If $Z_n \notin A_m$
then there must be $k \le n$ and $j \in \{ 1,\dots,d\}$ such that
$\hor(Z_{j,k}) < -m$. Thus, if we set
$$
S_n = \max \bigl\{ |\hor(Z_{j,k})| : j=1, \dots, d\,,\; k \le n \bigr\}\,,
$$
then $\Prob_o[Z_n \in A_m] \ge \Prob_o[S_n \le m]$. Since 
$\Hor(Z_n) = \bigl( \hor(Z_{1,n}), \dots, \hor(Z_{d,n}) \bigr)$
defines a translation invariant random walk on $\A_{d-1} \cong \Z^{d-1}$,
Lemma~1.2 of {\sc Alexopoulos~\cite{Al}} implies that there are constants
$C', D' > 0$ such that  
$$
\Prob_o[S_n \le m] \ge C' \exp(-D'n/m^2)\,.
$$
Therefore 
$$
p^{(2n)}(o,o) \ge \wt C \, 
\exp\Bigl( -\wt D \bigl( \tfrac{n}{m^2} + m + \log n \bigr)\Bigr)\,.
$$
Setting $m = \lfloor n^{1/3} \rfloor$, we obtain the claimed lower
estimate.
\end{proof} 

\medskip

\section{The Poisson boundary of random walk}\label{sect:Poisson}

In this section, we consider a larger class of random walks than those of
\S \ref{sect:spectrum}. On $\DL=\DL(q_1, \dots, q_d)$, consider a stochastic
transition matrix 
$$
P = \bigl( p(x,y) \bigr)_{x,y \in \DL}
$$
with the following properties:

\begin{enumerate}
\item[(i)] \emph{irreducibility:} For each $x, y \in \DL$ there is $N$ such
that the $(x,y)$-element of $P^N$ satisfies $p^{(N)}(x,y) > 0$.

\item[(ii)] \emph{group-invariance:} the group
$$
\Ga = \{ g \in \Af : p(gx,gy)=p(x,y) \;\text{for all}\; x,y \in \DL \}
$$
acts transitively on $\DL$, where $\Af$ is as in Proposition~\ref{pro:groupA}.  

\item[(iii)] \emph{finite first moment:} for the graph metric $d(\cdot,\cdot)$ of 
$\DL$, 
$$
m_1(P) = \sum_x p(o,x)d(o,x) < \infty\,.
$$
\end{enumerate}

Let $(Z_n)_{n \ge 0}$ be the random walk on $\DL$ governed by $P$, and let 
$\Prob_{x_0}$ be the probability measure on the associated trajectory space
for the initial point $Z_0=x_0$. Thus,
$$
\Prob_{x_0}[Z_{k+n}=y \mid Z_k = x] = p^{(n)}(x,y)
$$
whenever $\Prob_{x_0}[Z_k = x]= p^{(k)}(x_0,x) > 0$. 

\medskip

Since $\DL$ is a vertex-transitive graph with exponential growth,
the random walk is \emph{transient}, that is,
$$
\Prob_x[d(Z_n,Z_0) \to \infty]=1\,,
$$
see {\sc Woess~\cite[Theorem 5.13]{Wbook}}.
Our first question is whether one can give a more precise geometric description
of how $(Z_n)$ tends to infinity in $\DL$. 

\begin{dfn}\label{def:compactification} The 
\emph{geometric compactification} $\wh\DL$ of $\DL$ is the closure of $\DL$ in
$\prod_{j=1}^d \wh{T}_j$. The \emph{geometric boundary} of $\DL$ is
$$
\bd \DL= \wh \DL \setminus \DL = \left\{ \zeta = \zeta_1 \cdots \zeta_d \in 
\prod_{j=1}^d \wh{T}_j : \zeta_j = \om_j \;\text{for at least one}\;j\right\}.
$$
\end{dfn}

To justify the last fact, first observe that 
$\left(\prod_{j=1}^d \wh{T}_j\right) \setminus \left(\prod_{j=1}^d T_j\right)$ 
consists of all $\zeta = \zeta_1 \cdots \zeta_d \in \prod_{j=1}^d \wh{T}_j$ 
such that $\zeta_j \in \bd T_j$ for at least one $j$ (while it may
well be that some $\zeta_i$ are vertices in $T_i$).
Now, if such a $\zeta$ is the limit of a sequence 
$x_n=x_{1,n}\cdots x_{d,n} \in \DL$, then we can distinguish the
following cases: 

(1) If there are coordinates $i$ such that $|\hor(x_{i,n})|\to \infty$ then 
there must be some $j$ such that $\hor(x_{j,n})\to -\infty$ and consequently
$x_{j,n}\to \om_j$. 

(2) If $|\hor(x_{i,n})|$ is bounded for each $i$, then there still must
be some coordinate $j$ such that $d(o_j,x_{j,n}) \to \infty$. But in this case,
we also must have $\up(o_j,x_{j,n}) \to \infty$, see \eqref{updown}.
Therefore $x_{j,n}\to \om_j$. 

In general, up to a permutation of the coordinates, the situation will
be as follows: there are indices $1 \le r \le s \le d$ such that
$\zeta_j = \om_j$ for $1 \le j \le r$ (and $\hor(x_{j,n})$ may be bounded or
unbounded from below as well as from above~!), 
$\zeta_j \in \bd^*T_j$ for $r < j \le s$ (and $\hor(x_{j,n}) \to \infty$),
and $\zeta_j \in T_j$ for $s < j \le d$ (and $x_{j,n} = \zeta_j$ for all 
but finitely many $n$).

For a detailed description of boundary and convergence in the case
$d=2$ (two trees), see \cite[(5.3)]{BroWo} or {\sc Bertacchi~\cite{Ber}}.

For $j \in \{1,\dots,d\}$, let $\Ga_j$ be the image of $\Ga$ under
the projection $\Af \to \Aff(T_j)$. This group acts transitively on $T_j$.
Write $Z_{j,n}$ for the image of $Z_n$ under the projection $\DL \to T_j$.
Then $(Z_{j,n})$ is an irreducible random walk on $T_j$ whose transition
probabilities $p_j(x_j,y_j)$ are $\Ga_j$-invariant. We also consider
the image 
$$
\Hor(Z_n) = \bigl( \hor(Z_{1,n}), \dots, \hor(Z_{d,n}) \bigr) \in \A_{d-1}\,.
$$
This defines a random walk on $\A_{d-1}$ whose transition probabilities are
invariant under translation by elements of $\A_{d-1}$. In the same way,
$\bigl( \hor(Z_{j,n}) \bigr)_{n \ge 0}$ is a translation invariant
random walk on $\Z$. Its increments are i.i.d.\ integer random variables
with expected value
$$
\alpha_j = \alpha_j(P) = E_o\bigl(\hor(Z_{j,1})\bigr) 
= \sum_{x \in \DL} p(o,x) \hor(x_j)\,,
$$
which is finite, since $|\alpha_j| \le m_1(P)$. Note that 
$\sum_j \alpha_j = 0$, so that there must be $j$ with $\alpha_j \le 0$.
Applying the results of \cite{CaKaWo} and (for convergence to the boundary
when $\alpha_j=0$) {\sc Brofferio~\cite[Theorem 3.1]{Bro}}, one finds 
the following under the above assumptions (i)-(iii).

\begin{pro}\label{pro:convergence}
For each $j \in \{1, \dots, d\}$, we have for the distance on $T_j$
$$
\lim_{n \to \infty} \frac{d(Z_{j,n},Z_{0,n})}{n} = \alpha_j \quad
\Prob_x\text{-almost surely}\,,
$$
and there is a $\bd T_j$-valued random variable $Z_{j,\infty}$
such that
$$
\lim_{n \to \infty} Z_{j,n} = Z_{j,\infty} \quad
\Prob_x\text{-almost surely}
$$
in the topology of $\wh{T}$, for every $x \in \DL$. One has the 
following.\\[5pt]
{\rm (a)} If $\alpha_j \le 0$ then $Z_{j,\infty} = \om_j$ almost surely.\\[5pt]
{\rm (b)} If $\alpha_j > 0$ then $Z_{j,\infty}$ is a $\bd^*T_j$-valued
random variable whose $\Prob_x$-distribution is a continuous measure
(i.e., it carries no point mass)
supported by the whole of $\bd^*T_j$.
\end{pro}
  
\begin{proof} The results in \cite{CaKaWo} and \cite{Bro} are formulated
for random walks on $\Aff(\T_q)$, that is, for $R_no$, where 
$R_n=g_0 X_1 \cdots X_n$
(group product) and  the $X_n$ are i.i.d.\ group-valued random variables. 
The only point that we have to clarify is that our situation arises as
a special case of that form.

There is a standard method for turning a $\Ga$-invariant random walk
on a transitive graph into $R_no$, where $(R_n)$ is a random walk on
$\Ga$; see e.g. {\sc Kaimanovich and Woess~\cite[Prop. 2.15]{KaWo}}. 
Indeed, if $dg$ is the left Haar measure on the 
(locally compact, totally disconnected) group $\Ga$, then
\begin{equation}\label{eq:mu}
\mu(dg) = p(o,go)\,dg
\end{equation} 
defines a probability measure on $\Ga$. Now let
$(X_n)_{n\ge 1}$ be i.i.d. $\mu$-distributed in $\Ga$. Then 
$Z_n=g_0 X_1 \cdots X_n o$ defines a Markov chain on the graph with
$Z_0=g_0o$ and transition matrix~$P$. 
\end{proof}

Let us write $\Lambda_j(\alpha_j) = \{ \om_j \}$ if $\alpha_j \le 0$,
and $\Lambda_j(\alpha_j) = \bd^*T_j$ if $\alpha_j > 0$, and define the
following Borel subset of $\bd\DL$.
\begin{equation}\label{eq:Lambda} 
\Lambda = \prod_{j=1}^d \Lambda_j(\alpha_j).
\end{equation} 

\begin{cor}\label{cor:bdry} Under assumptions (i)-(iii), there is
a $\Lambda$-valued random variable $Z_{\infty}$ such that 
$$
\lim_{n \to \infty} Z_{n} = Z_{\infty} \quad
\Prob_x\text{-almost surely}
$$
in the topology of $\wh\DL$, for every $x \in \DL$. 
The $\Prob_x$-distribution $\nu_x$ of $Z_{\infty}$ is a
Borel probability measure supported by the whole of $\Lambda$.
It is a continuous measure unless $\Lambda = \{ \om_1\cdots \om_d \}$
is trivial, which happens precisely when $\alpha_j=0$ for all $j$.

The measures $\nu_x$ satisfy $\nu_x = \sum_y p(x,y)\nu_y$ and are
mutually absolutely continuous. 
\end{cor} 

(The last statement is immediate by factoring through the first step
of the random walk, and using irreducibility.)

Thus, we see that the space $\Lambda$ together with the family of \emph{harmonic
measures} $\nu_{\bullet}=(\nu_x)_{x \in \DL}$ is a \emph{boundary} of the random walk
$(Z_n)$ in the sense of \cite{KaWo}, a notion going back to
{\sc Furstenberg~\cite{Fu}}. Indeed, the group $\Ga$ defined in (ii)
acts on $\Lambda$, i.e., $\Lambda$ is a $\Ga$-space, and we have the convolution
identity $\mu*\nu_o = \nu_o$ for the probability measure $\mu$ on $\Ga$
defined in \eqref{eq:mu}. This identity holds since 
$\de_g * \nu_o = \nu_{go}$ for every $g \in \Ga$.

We now want to decide whether this is the ``best'' (biggest) model, as a
measure space, for distinguishing limit points at infinity of our
random walk. To formulate this question more precisely, we recall that
a $P$-harmonic function is a function $h: \DL \to \R$ such that $h(x)
= \sum_y p(x,y)h(y)$. Now, if $\varphi \in
L^{\infty}(\Lambda,\nu_{\bullet})$, then $$
h(x) = \int_{\Lambda}
\varphi\,d\nu_x $$
defines a mapping from
$L^{\infty}(\Lambda,\nu_{\bullet})$ to the space of all bounded
harmonic functions. When the mapping is bijective,
$(\Lambda,\nu_{\bullet})$ is called the \emph{Poisson boundary} of the
random walk. The latter is unique up to isomorphism of measure spaces.
For precise details in the present setting of a group-invariant random
walk on a transitive graph, we refer to~\cite{KaWo}; additional
information can be found in the references given there. The following
theorem should be compared with the main result (based on a different
method) of {\sc Brofferio~\cite{Bro2}}, who considers the Poisson
boundary of random walks on finitely generated groups of affine
mappings with rational coefficients: those groups act on a
``weighted'' horocyclic product of distinct trees with hyperbolic
upper half plane.

\begin{thm}\label{thm:poisson} Under assumptions (i)-(iii), 
$(\Lambda,\nu_{\bullet})$ is the Poisson boundary of  
the random walk $(Z_n)$ on $\DL$.
\end{thm}

\begin{proof} For the proof, we apply the \emph{ray criterion} of
{\sc Kaimanovich}, see \cite[Thm. 5.18]{KaWo}. Formulated in our specific
terms of $\Lambda$ and $\DL$, it says the following.

\begin{enumerate}
\item[]
Suppose that $\Lambda$ is a $\Ga$-space carrying a family of probability
measures $\nu_{\bullet} = (\nu_x)_{x \in \DL}$ such that 
$\de_g*\nu_o = \nu_{go}$ and $\nu_o = \sum_x p(o,x)\nu_x$.
Then $(\Lambda, \nu_{\bullet})$ is the Poisson boundary of the random
walk $(Z_n)$, provided that there is a sequence of measurable 
mappings $\Pi_n: \Lambda \to \DL$ such that, in the graph metric of
$\DL$, 
$$
\lim_{n\to\infty} \frac{1}{n}d\bigl(Z_n, \Pi_n(Z_{\infty})\bigr) = 0
\quad \Prob_o\text{-almost surely.}
$$
\end{enumerate}

We now construct $\Pi_n$.  We know that $\alpha_j \le 0$ for some $j$,
and assume without loss of generality that $\alpha_d \le 0$.  Let now
$\xi = \xi_1 \cdots \xi_d \in \Lambda$ be given, so that $\xi_j \in
\bd T_j$ for each $j$, and $\xi_j = \om_j$ for at least one $j$.  By
Proposition \ref{pro:convergence}, $\xi_d = \om_d$.

For integer $k \ge 0$, denote by $\xi_j(k)$ the element on the geodesic 
$\geo{o_j\,\xi_j}$ at distance $k$ from $o_j$. 
Similarly, if $\xi_j \in \bd^*T_j$ and $k \ge 0$,
denote by $\xi_j[k]$ the element on the geodesic 
$\geo{o_j\,\xi_j}$ which lies on the horocycle $H_k^j$ of $T_j$ and, in case 
there are two such elements (only when $k=0$), is furthest possible from 
$o_j$.  Finally, if $\xi_j=\om_j$ and $k \le 0$ then the geodesic
$\geo{o_j\,\xi_j}$ intersects $H_k$ at the point $\xi_j[k]=\om_j(|k|)$.
See Figure~8. 
$$
\beginpicture 

\setcoordinatesystem units <6mm,8mm>

\setplotarea x from -6 to 6, y from -2 to 6

\arrow <6pt> [.2,.67] from 0 0 to 6 6
\arrow <6pt> [.2,.67] from 2 2 to 6 -2

\setdots <3pt>

\putrule from -4.6 3 to 8 3
\putrule from -4.6 0 to 8 0
\putrule from -4.6 -1 to 8 -1

\put {$H_{k}^j\; (k > 0)$} [r] at -4.7 -1
\put {$H_{m}^j\; (m < 0)$} [r] at -4.7  3
\put {$H_0^j$} [r] at -4.7 0

\put {$o_j$} [rb] at -0.1 0.1 
\put {$\om_j$} [lb] at 6.1 6.1
\put {$\xi_j$} [lt] at 6.1 -2.1

\put {$\xi_j(k)$} [rb] at 0.95 1.15
\put {$\xi_j[0]$} [lb] at 4.05 0.15
\put {$\xi_j[k]$} [lb] at 5.05 -0.85

\put {$\om_j[m]=\om_j(|m|)$} [rb] at 2.95 3.15

\multiput {\scriptsize $\bullet$} at 0 0  1 1  3 3
                                       4 0  5 -1 /

\endpicture
$$

\vspace{-.3cm}

\begin{center}
\centerline\emph{Figure 8}
\end{center}

\vspace{.1cm}

\noindent Note that 
\begin{equation}\label{eq:geobound}
d\bigl(\xi_j[k],\xi_j(|k|)\bigr) = d(o_j,\xi_j[0])\,.
\end{equation}
We define 
$$
\begin{gathered}
\Pi_n(\xi) = x_1 \cdots x_d\,, \quad \text{where} \quad
x_j = \Pi_{j,n}(\xi_j) = \xi_j[k_j] \quad \text{with} \\ 
k_j=\lceil \alpha_j n \rceil\,,\; j=1,\dots,d-1\,,\AND
k_d = -k_1- \ldots - k_{d-1}\,.
\end{gathered}
$$
Note that this is well defined: when $\alpha_j > 0$ we have $\xi_j
\in \bd^*T_j$ and $k_j > 0$, so that $H_{k_j}^j$ intersects
$\geo{o_j\,\xi_j}$ in a unique point. When $j < d$ and $\alpha_j \le
0$ we have $\xi_j = \om_j$ and $k_j\le 0$, so that $H_{k_j}^j$ also
intersects $\geo{o_j\,\om_j}$ in a unique point. This also holds for
$j=d$, since $k_d \le -\alpha_1 n -\ldots - \alpha_{d-1}n = \alpha_dn
\le 0$.

We now apply the results and methods of \cite{CaKaWo}.
Since $m_1(P) < \infty$, we also have $m_1(P_j) < \infty$, where $P_j$ is the
transition matrix of $(Z_{j,n})_{n \ge 0}$ on $T_j$. This implies via the law of large
numbers that
$$
\frac{1}{n} d(Z_{j,n+1},Z_{j,n}) \to 0 \AND
\frac{1}{n} \hor(Z_{j,n}) \to \alpha_j \quad \Prob_o\text{-almost surely.}
$$
Applying \cite[Def.~1 and Prop.~1]{CaKaWo}, one obtains that $(Z_{j,n})_n$ is
almost surely a \emph{regular sequence} in $T_j$, that is,
$$
\frac{1}{n} d\bigl(Z_{j,n},Z_{j,\infty}(\lceil \alpha_j n \rceil)\bigr) 
 \to 0\quad \Prob_o\text{-almost surely.} 
$$
Now \eqref{eq:geobound} implies that in the metric of $T_j$ also
$$
\frac{1}{n} d\bigl(Z_{j,n},\Pi_{j,n}(Z_{j,\infty})\bigr) 
 \to 0\quad \Prob_o\text{-almost surely.} 
$$
Since for arbitrary $x, y \in \DL$ we have 
$d(x,y) \le \sum_{j=1}^d d(x_j,y_j)$, we obtain that in the metric of $\DL$
$$
\frac{1}{n} d\bigl(Z_n,\Pi_n(Z_{\infty})\bigr) 
 \to 0\quad \Prob_o\text{-almost surely,} 
$$
so that the  criterion applies.
\end{proof}

Since $Z_{j,n}\to \om_j$ almost surely, when $\alpha_j \le 0$, we may 
omit this part in the description of the Poisson boundary, and
can consider the harmonic measures $\nu_x$ as measures on the product
of all $\bd^*T_j$ where $\alpha_j > 0$.

\begin{cor}\label{cor:bdproduct}
The Poisson boundary of the random walk $(Z_n)$ satisfying (i)-(iii)
can be identified with the space 
$$
\prod_{j: \alpha_j > 0} \bd^*T_j\,.
$$
When $\alpha_j=0$ for all $j$ then the Poisson boundary is trivial, and
all bounded harmonic functions are constant.

In particular, for SRW on $\DL(q_1,\dots,q_d)$ we have, setting 
$\bar q=(q_1 + \dots + q_d)/d$, that the Poisson boundary can be identified with
$$
\prod_{j: q_j > \bar q} \bd^*T_j\,;
$$
it is trivial precisely when all $q_j$ coincide.
\end{cor} 

Indeed, for SRW, we compute $\alpha_j = (q_j - \bar q)d/D$, where 
$D$ is the vertex degree in $\DL$.

\vspace{.6cm}

\textbf{Acknowledgments.} We acknowledge an email exchange
with Brian J. McCartin on the spectrum of the Laplacian on equilateral
triangles. We are particularly indebted to Christophe Pittet for
an extremely useful discussion about the construction of co-compact
lattices in the isometry group of $\DL_3(q)$. Fruitful conversations
with Andr\'e Henriques and Christian W\"uthrich are gratefully
acknowledged.


\end{document}